\documentclass[10pt]{article}

\makeatletter
\newfont{\footsc}{cmcsc10 at 8truept}
\newfont{\footbf}{cmbx10 at 8truept}
\newfont{\footrm}{cmr10 at 10truept}
\renewcommand{\ps@plain}{%
\renewcommand{\@oddfoot}{\footsc the electronic journal of combinatorics
  {\footbf 12} (2005), \#R00\hfil\footrm\thepage}}
\makeatother
\title{A Gessel--Viennot-Type Method for Cycle Systems\\ in a Directed Graph}

\author{Christopher R. H. Hanusa\\
\small Department of Mathematical Sciences\\[-0.8ex]
\small Binghamton University, Binghamton, New York, USA\\[-0.8ex]
\small \texttt{chanusa@math.binghamton.edu}}

\date{\small 
Submitted: Nov 28, 2005;  Accepted: Jan 2, 2005; Published: Jan 3, 2005\\
\small Mathematics Subject Classifications: Primary 05B45, 05C30;\\ 
\small Secondary 05A15, 05B20, 05C38, 05C50, 05C70, 11A51, 11B83, 15A15, 15A36, 52C20 \\
\small Keywords: directed graph, cycle system, path system, walk system, Aztec diamond, \\ Aztec pillow, 
Hamburger Theorem, Kasteleyn--Percus, Gessel--Viennot, Schr\"oder numbers}

\usepackage{amsthm}
\usepackage{epsfig}
\usepackage{amssymb}
\usepackage{amsmath}

\newtheorem{theorem}{Theorem}[section]
\newtheorem{corollary}[theorem]{Corollary}
\newtheorem{conjecture}[theorem]{Conjecture}
\newtheorem{lemma}[theorem]{Lemma}
\theoremstyle{definition}
\newtheorem{definition}[theorem]{Definition}

\theoremstyle{remark}

\def\CC{{\mathcal{C}}} \def\C{{\mathcal{W}}} \def\F{{\mathcal{F}}} \def\O{{\mathcal{O}}} \def\P{{\mathcal{P}}}
\def\i{{\iota}}  \def\phi{\varphi} \def\ep{\varepsilon} \def\lam{\lambda}  
\def\ell{{l}}   \def\wt{\text{wt}} \def\sgn{\text{sgn}}
\def\row{\rightarrow}  \def\lrow{\longrightarrow}

\begin{document}
\maketitle

%%%%%%%%%%%%%%%%%%%%%%%%%%%%%%%%%%%%%%%%%%%%%%%%%%%%%%%%%%
\begin{abstract}
We introduce a new determinantal method to count cycle systems in a directed graph that generalizes Gessel 
and Viennot's determinantal method on path systems.  The method gives new insight into 
the enumeration of domino tilings of Aztec diamonds, Aztec pillows, and related regions.  
\end{abstract}

%%%%%%%%%%%%%%%%%%%%%%%%%%%%%%%%%%%%%%%%%%%%%%%%%%%%%%%%%%%%%%%%
\section{Introduction}

\begin{figure}
\begin{center}
  \epsfig{figure=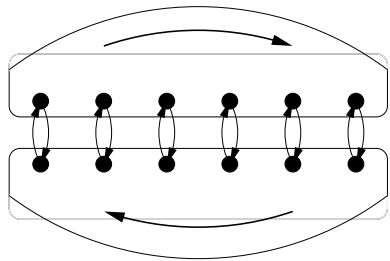}
  \caption{A hamburger graph}
  \label{Hamburger}
\end{center}
\end{figure}

In this article, we present an analogue of the Gessel--Viennot method for counting cycle systems on a type of 
directed graph we call a hamburger graph.  A {\em hamburger graph} $H$ is made up of two acyclic graphs $G_1$ and 
$G_2$ and a connecting edge set $E_3$ with the following properties.  The graph $G_1$ has $k$ distinguished 
vertices $\{v_1,\hdots,v_k\}$ with directed paths from $v_i$ to $v_j$ only if $i< j$. The graph $G_2$ has $k$ 
distinguished vertices $\{w_{k+1},\hdots,w_{2k}\}$ with directed paths from $w_i$ to $w_j$ only if $i> j$.  The 
edge set $E_3$ connects each vertex $v_i$ to vertex $w_{k+i}$ and vice versa.  (See Figure \ref{Hamburger} for a 
visualization.)  Hamburger graphs arise naturally in the study of Aztec diamonds, as explained in Section 
\ref{Apps}.

%\begin{figure}[b]

The Gessel--Viennot method is a determinantal method to count path systems in an acyclic directed graph $G$ with 
$k$ sources $s_1,\hdots,s_k$ and $k$ sinks $t_1,\hdots,t_k$.  A {\em path system} $\P$ is a collection of $k$ 
vertex-disjoint paths, each one directed from $s_i$ to $t_{\sigma(i)},$ for some permutation $\sigma\in S_k$ 
(where $S_k$ is the symmetric group on $k$ elements).  Call a path system $\P$ {\em positive} if the sign of this 
permutation $\sigma$ satisfies $\sgn(\sigma)=+1$ and {\em negative} if $\sgn(\sigma)=-1$.  Let $p^+$ be the number 
of positive path systems and $p^-$ be the number of negative path systems.

Corresponding to this graph $G$ is a $k\times k$ matrix $A=(a_{ij})$, where $a_{ij}$ is the number of paths from 
$s_i$ to $t_j$ in $G$.  
The result of Gessel and Viennot states that $\det A = p^+-p^-$.  
The Gessel--Viennot method was introduced in \cite{GV1,GV2}, and has its roots in works by Karlin and McGregor 
\cite{KarMcG} and Lindstr\"om \cite{Lindstrom}.  A nice exposition of the method and applications is given in the 
article by Aigner \cite{Aigner}.  

%For additional applications of the Gessel--Viennot method, see \cite{Destain,Mays}.

This article concerns a similar determinantal method for counting cycle systems in a hamburger graph $H$.  A 
{\em cycle system} $\CC$ is a collection of vertex-disjoint directed cycles in $H$.  Let $\ell$ be the number of 
edges in $\CC$ that travel from $G_2$ to $G_1$ and let $m$ be the number of cycles in $\CC$.  Call a cycle system 
{\em positive} if $(-1)^{\ell+m}=+1$ and {\em negative} if $(-1)^{\ell+m}=-1$.  Let $c^+$ be the number of 
positive cycle systems and $c^-$ be the number of negative cycle systems.
%
%\medskip
Corresponding to each hamburger graph $H$ is a $2k\times 2k$ block matrix $M_H$ of the form
\begin{equation*}
M_H=
\left[
\begin{array}{cc}
A & I_k \\
-I_k & B  \\
\end{array}
\right],
\end{equation*}
where in the upper triangular matrix $A=(a_{ij})$, $a_{ij}$ is the number of paths from $v_i$ to $v_j$ in 
$G_1$ and in the lower triangular matrix $B=(b_{ij})$, $b_{ij}$ is the number of paths
from $w_{k+i}$ to $w_{k+j}$ in $G_2$.  This matrix $M_H$ is referred to as a {\em hamburger matrix}.

{\theorem[The Hamburger Theorem] If $H$ is a hamburger graph, then $\det M_H = c^+-c^-$. \label{HT}}

\medskip
A hamburger graph $H$ is called {\em strongly planar} if there is a planar embedding of $H$ that sends $v_i$ to 
$(i,1)$ and $w_{k+i}$ to $(i,-1)$ for all $1\leq i\leq k$, and keeps edges of $E_1$ in the half-space $y\geq 1$ 
and edges of $E_2$ in the half-space $y\leq -1$.  This definition suggests that $G_1$ and $G_2$ are ``relatively'' 
planar in $H$, a stronger condition than planarity of $H$.  Notice that when $H$ is strongly planar, each cycle 
must use 
exactly one edge from $G_2$ to $G_1$.  Hence, the sign of every cycle system is $+1$. This implies the following 
corollary.

{\corollary If $H$ is a strongly planar hamburger graph, $\det M_H = c^+$.}

\begin{figure}
\begin{center}  
  \epsfig{figure=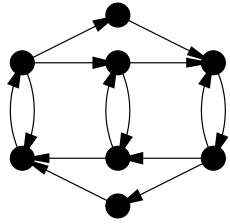} 
\end{center}
  \caption{A simple hamburger graph $H$}
  \label{NiceGraph}
\end{figure}

\medskip 
The following simple example serves to guide us.  Consider the two graphs $G_1=(V_1,E_1)$ and $G_2=(V_2,E_2)$, 
where $V_1=\{v_1,v_2,v_3\}$, $V_2=\{w_4,w_5,w_6\}$, $E_1=\{v_1\row v_2, v_2\row v_3, v_1\row v_3\}$, and 
$E_2=\{w_6\row w_5, w_5\row w_4, w_6\row w_4\}$.  Our hamburger graph $H$ will be the union of $G_1$, 
$G_2$, and the edge set $E_3$ consisting of edges $e_i:v_i\row w_{k+i}$ and $e_i':w_{k+i}\row v_i$.  In this 
example, $k=3$ and $H$ is strongly planar.  Figure \ref{NiceGraph} gives a graphical representation of $H$.

\begin{figure}
\begin{center}  
  \epsfig{figure=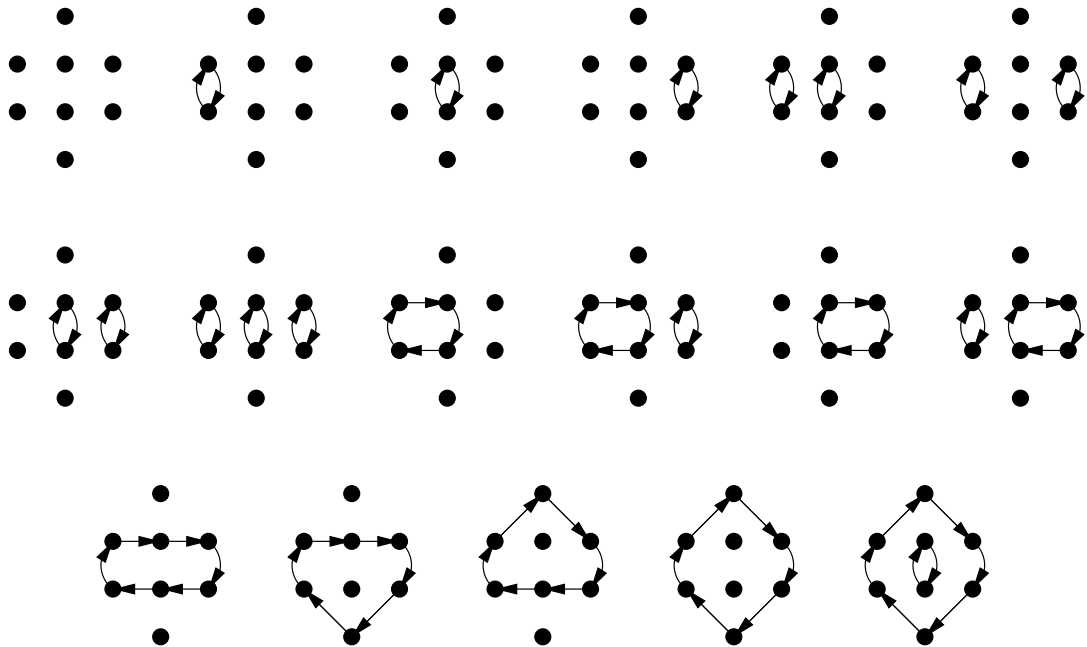} 
\end{center}
  \caption{The seventeen cycle systems for the hamburger graph in Figure \ref{NiceGraph}}
  \label{SysCycles}
\end{figure}

In this example, the hamburger matrix $M_H$ equals
\begin{equation*}
M_H=
\left[
\begin{array}{cccccc}
 1 & 1 & 2 & 1 & 0 & 0 \\
 0 & 1 & 1 & 0 & 1 & 0 \\
 0 & 0 & 1 & 0 & 0 & 1 \\
-1 & 0 & 0 & 1 & 0 & 0 \\
0 & -1 & 0 & 1 & 1 & 0 \\
0 & 0 & -1 & 2 & 1 & 1 \\
\end{array}
\right].
\end{equation*}
The determinant of $M_H$ is 17, corresponding to the seventeen cycle systems (each with sign +1) in Figure \ref{SysCycles}.  

\medskip

The graph that inspired the definition of a hamburger graph comes from the work of Brualdi and Kirkland 
\cite{BruKirk}, in which they give a new proof that the number of domino tilings of the Aztec diamond is 
$2^{n(n+1)/2}$. An {\em Aztec diamond}, denoted by $AD_n$, is the union of the $2n(n+1)$ unit squares with 
integral vertices $(x,y)$ such that $|x|+|y|\leq n+1$.  An {\em Aztec pillow}, as it was initially presented in 
\cite{Propp}, is also a rotationally symmetric region in the plane.  On the top left boundary, however, the steps 
are composed of three squares to the right for every square up.  Another definition is that Aztec pillows are the 
union of the unit squares with integral vertices $(x,y)$ such that $|x+y|<n+1$ and $|3y-x|<n+3$.
As with Aztec diamonds, we denote the Aztec 
pillow with $2n$ squares in each of the central rows by $AP_{n}$.  In Section \ref{Generalizations}, we extend the 
notion of Aztec pillows having steps of length 3 to ``odd pillows''---those that have steps that are of a constant 
odd length.  The integral vertices $(x,y)$ of the unit squares in $q$-pillows for $q$ odd satisfy $|x+y|<n+1$ and 
$|qy-x|<n+q$.

\begin{figure}
\begin{center}
\begin{tabular}{ccc}
  \epsfig{figure=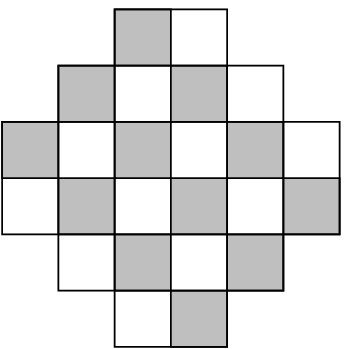} &
  \epsfig{figure=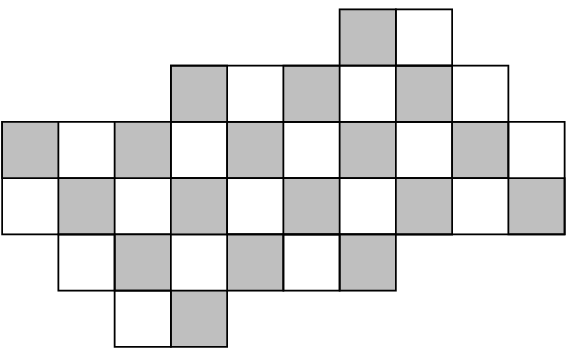} &
  \epsfig{figure=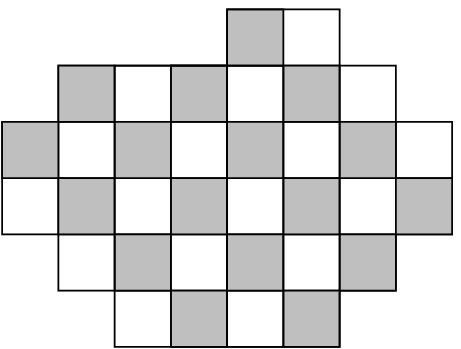} 
\end{tabular}
  \caption{Examples of an Aztec diamond, an Aztec pillow, and a generalized Aztec pillow}
  \label{AztecEx}
\end{center}
\end{figure}

We introduce the idea of a {\em generalized Aztec pillow}, where the steps on all diagonals are of possibly 
different odd lengths.  More specifically, a generalized Aztec pillow is a horizontally convex and vertically 
convex region such that the steps both up and down in each diagonal have an odd number of squares horizontally for 
every one square vertically.  A key fact that we will use is that any generalized Aztec pillow can be recovered 
from a large enough Aztec diamond by the placement of horizontal dominoes.  See Figure \ref{AztecEx} for examples 
of an Aztec diamond, an Aztec pillow, and a generalized Aztec pillow.

Brualdi and Kirkland prove the formula for the number of domino tilings of an Aztec diamond by creating an 
associated digraph and counting its cycle systems, manipulating the digraph's associated Kasteleyn--Percus matrix 
of order $n(n+1)$.  To learn about Kasteleyn theory and Kasteleyn--Percus matrices, start with Kasteleyn's 1961 
work \cite{Kasteleyn} and Percus's 1963 work \cite{Percus}.  The Hamburger Theorem proves that we can count the number 
of domino tilings of an Aztec diamond with a much smaller determinant, of order $2n$.  An analogous reduction in 
determinant size (from order $O(n^2)$ to order $O(n)$) occurs for all regions to which this theorem applies, 
including generalized Aztec pillows.  In addition, whereas Kasteleyn theory applies only to planar graphs, there 
is no planarity restriction for hamburger graphs.  For this reason, the Hamburger Theorem gives a new counting 
method for cycle systems in some non-planar graphs.

In Section \ref{Outline}, we present an overview of the proof of the Hamburger Theorem, including the key lemmas 
involved.  The necessary machinery is built up in Section \ref{Additional} to complete the proof in Section 
\ref{HamProof}.  Section \ref{Apps} presents applications of the Hamburger Theorem to Aztec diamonds, Aztec 
pillows, and generalized Aztec pillows.  Section \ref{Generalizations} concludes with a counterexample to the most 
natural generalization of the Hamburger Theorem and an extension of Propp's Conjecture on Aztec pillows.

%%%%%%%%%%%%%%%%%%%%%%%%%%%%%%%%%%%%%%%%%%%%%%%%%%%%%%%%%%%%%%%%
\section{Outline of the Proof of the Hamburger Theorem}
\label{Outline}
%%%%%%%%%%%%%%%%%%%%%%%%%%%%%%%%%%%%%%%%%%%%%%%%%%%%%%%%%%%%%%%%

\subsection{The Hamburger Theorem}

Like the proof of the Gessel--Viennot method, the proof of the Hamburger Theorem hinges on cancellation of terms 
in the permutation expansion of the determinant of $M_H$.  In the proof, we must allow closed directed walks in 
addition to cycles.  We must also allow {\em walk systems}, arbitrary collections of closed directed walks, since 
they can and will appear in the permutation expansion of the hamburger determinant.  We call a walk system {\em 
simple} if the set of walks visits no vertex more than once.  We call a cycle of the form $c:v_i \row w_{i+k} \row 
v_i$ a {\em 2-cycle}.

Each signed term in the permutation expansion of the hamburger determinant is the contribution of many signed walk 
systems $\C$.  Walk systems that are not cycle systems will all cancel out in the determinant expansion.  We will 
show this in two steps. We start by considering walk systems that are not simple.  If this is the case, one of the 
two following properties MAY hold.

\medskip
\noindent
{\bf Property 1.} The walk system contains a walk that has a self-intersection.  

\medskip
\noindent
{\bf Property 2.} The walk system has two intersecting walks, neither of which is a 2-cycle.  

\medskip
The following lemma shows that the contributions of walk systems satisfying either of these two properties cancel 
in the permutation expansion of the determinant of $M_H$.

{\lemma The set of all walk systems $\C$ that satisfy either Property 1 or Property 2 can be 
partitioned into equivalence classes, each of which contributes a net zero to the permutation expansion of the 
determinant of $M_H$. \label{Lem1}}

\medskip
The proof of Lemma \ref{Lem1} uses a generalized involution principle.  Walk systems cancel in families based on 
the their ``first'' intersection point.  

\medskip
The remainder of the cancellation in the determinant expansion is based on the concept of a minimal walk system; 
we motivate this definition by asking the following questions.  What kind of walk 
systems does the permutation expansion of the hamburger determinant generate, and how is this different from our 
original notion of cycle systems that we wanted to count in the introduction?  The key difference is that the same 
collection of walks can be generated by multiple terms in the determinantal expansion of $M_H$; whereas, we 
would only want to count it once as a cycle system.  This redundancy arises when the walk visits three 
distinguished vertices in $G_1$ without passing via $G_2$ or vice versa.  We illustrate this notion with the 
following example.

Consider the second cycle system in the third row of Figure \ref{SysCycles}, consisting of one solitary directed 
cycle.  Since this cycle visits vertices $v_1$, $v_2$, $v_3$, $w_6$, and $w_4$ in that order, it contributes a 
non-zero weight in the permutation expansion of the determinant corresponding to the term $(12364)$ in $S_6$.  
Notice that this cycle {\em also} contributes a non-zero weight in the permutation expansion of the determinant 
corresponding to the term $(1364)$.  We see this since our cycle follows a path from $v_1$ to $v_3$ (by way of 
$v_2$), returning to $v_1$ via $w_6$ and $w_4$.  We must deal with this ambiguity.  We introduce the idea of a 
{\em minimal} permutation cycle, one which does not include more than two successive entries with values between 
$1$ and $k$ or between $k+1$ and $2k$.  We see that $(1364)$ is minimal while $(12364)$ is not.

We notice that walk systems arise from permutations, so it is natural to think of a walk system first as a 
permutation, and afterward a collection of walks that ``follow'' the permutation.  This is the idea of a walk 
system--permutation pair (or WSP-pair for short) that is presented in Section \ref{WSPPair}.  From the idea of a 
minimal permutation cycle, we define a {\em minimal walk} to have as its base permutation a minimal permutation 
cycle, and a {\em minimal walk system} to be composed of only minimal walks. Since our original goal was to count 
``cycle systems'' in a directed graph, we realize we need to be precise and instead count ``simple minimal walk 
systems''.  This leads to the second part of the proof of the Hamburger Theorem.

Given a walk system that satisfies neither Property 1 nor Property 2 and is neither simple nor minimal, at least 
one of the two following properties MUST hold.  

\medskip
\noindent
{\bf Property 3.} The walk system has two intersecting walks, one of which is a 2-cycle.  

\medskip
\noindent
{\bf Property 4.} The walk system is not minimal.  

\medskip
The following lemma shows that the contributions of walk systems satisfying either of these 
new properties cancel in the permutation expansion of the determinant of $M_H$.

{\lemma The set of all walk systems $\C$ that satisfy neither Property 1 nor Property 2 and that satisfy either 
Property 3 or Property 4 can be partitioned into equivalence classes, each of which contributes a net zero to the 
permutation expansion of the determinant of $M_H$. \label{Lem2}}

\medskip
The proof of Lemma \ref{Lem2} is also based on involutions.  Walk systems cancel in families built from an index 
set containing the set of all 2-cycle intersections and non-minimalities.

\medskip
If a walk system satisfies none of the conditions of Properties 1 through 4, then it is indeed a simple minimal 
walk system, or in other words, a cycle system.  The cancellation from the above sets of families gives that only 
cycle systems contribute to the permutation expansion of the determinant of $M_H$.  This contribution is the 
signed weight of each cycle system, so the determinant of $M_H$ exactly equals $c^+-c^-$.  
Theorem \ref{HT} follows from Lemmas \ref{Lem1} and \ref{Lem2} in Section \ref{HamProof}. \hfill$\bullet$

\subsection{The Weighted Hamburger Theorem}
\label{Sec:WHT}

There is also a weighted version of the Hamburger Theorem, and it will be under this generalization that Lemmas 
\ref{Lem1} and \ref{Lem2} are proved.  We allow {\em weights} $\wt(e)$ on the edges of the hamburger graph; the 
simplest weighting, which counts the number of cycle systems, assigns $\wt(e)\equiv 1$.  We require that 
$\wt(e_i)\wt(e_i')=1$ for all $2\leq i \leq k-1$, but we do not require this condition for $i=1$ nor for $i=k$.  
Define the $2k\times 2k$ {\em weighted hamburger matrix} $M_H$ to be the block matrix
\begin{equation}
M_H=
\left[
\begin{array}{cc}
A & D_{1} \\
-D_{2} & B  \\
\end{array}
\right].
\label{BlockMatrix}
\end{equation}
In the upper-triangular $k\times k$ matrix $A=(a_{ij})$, $a_{ij}$ is the sum of the products of the weights of 
edges over all paths from $v_i$ to $v_j$ in $G_1$.  In the lower triangular $k\times k$ matrix $B=(b_{ij})$, 
$b_{ij}$ is the sum of the products of the weights of edges over all paths from $w_{k+i}$ to $w_{k+j}$ in $G_2$.  
The diagonal $k\times k$ matrix $D_{1}$ has as its entries $d_{ii}=\wt(e_i)$ and the diagonal $k\times k$ matrix 
$D_{2}$ has as its entries $d_{ii}=\wt(e_i')$.  Note that when the weights of the edges in $E_3$ are all 1, these 
matrices satisfy $D_{1}=D_{2}=I_k$.

We wish to count vertex-disjoint unions of weighted cycles in $H$.  
In any hamburger graph $H$, there are two possible types of cycle.  There are $k$ {\em $2$-cycles} 
\begin{equation*}
c:v_i\stackrel{e_i}{\lrow} w_{k+i} \stackrel{e_i'}{\lrow} v_i
\end{equation*}
and many more {\em general cycles} that alternate between $G_1$ and $G_2$.  We can think of a general cycle
as a path $P_1$ in $G_1$ connected by an edge $e_{1,1}\in E_3$ to a path $Q_1$ in $G_2$, which in turn
connects to a path $P_2$ in $G_1$ by an edge $e_{1,2}'$, continuing in this fashion until arriving at a final path
$Q_\ell$ in $G_2$ whose terminal vertex is adjacent to the initial vertex of $P_1$.  We write
\begin{equation*}
c: P_1 \stackrel{e_{1,1}}{\lrow} 
Q_1 \stackrel{e_{1,2}'}{\lrow} 
P_2 \stackrel{e_{2,1}}{\lrow} \cdots 
\stackrel{e_{\ell, 1}}{\lrow} 
P_{\ell} \stackrel{e_{\ell, 2}'}{\lrow} Q_{\ell}.
\end{equation*}
For each cycle $c$, we define the weight $\wt(c)$ of $c$ to be the product of the weights of all edges traversed 
by $c$: 
\begin{equation*}
\wt(c)=\prod_{e\in c} \wt(e).
\end{equation*}

We define a {\em weighted cycle system} to be a collection $\CC$ of $m$ vertex-disjoint cycles.  We again define 
the {\em sign} of a weighted cycle system to be $\sgn(\CC)= (-1)^{\ell+m}$, where $\ell$ is the total number of 
edges from $G_2$ to $G_1$ in $\CC$.  We say that a weighted cycle system $\CC$ is {\em positive} if $\sgn(\CC)=+1$ 
and {\em negative} if $\sgn(\CC)=-1$.  For a hamburger graph $H$, let $c^+$ be the sum of the weights of positive 
weighted cycle systems, and let $c^-$ be the sum of the weights of negative weighted cycle systems.

\noindent
{\theorem[The weighted Hamburger Theorem] The determinant of the weighted hamburger matrix $M_H$ equals $c^+-c^-$. \label{WHT}}

\medskip
As above, Theorem \ref{WHT} follows from Lemmas \ref{Lem1} and \ref{Lem2}.  The proofs will be presented after 
developing the following necessary machinery.  

%%%%%%%%%%%%%%%%%%%%%%%%%%%%%%%%%%%%%%%%%%%%%%%%%%%%%%%%%%%%%%%%
\section{Additional Definitions}
\label{Additional}

%\medskip
\subsection{Edge Cycles and Permutation Cycles}

In the proof of the Hamburger Theorem, there are two distinct mathematical objects that have the name ``cycle''.  
We have already mentioned the type of cycle that appears in graph theory.  There, a (simple) cycle in a directed 
graph is a closed directed path with no repeated vertices.

Secondly, there is a notion of cycle when we talk about permutations.  If $\sigma\in S_n$ is a permutation, we can 
write $\sigma$ as the product of disjoint cycles $\sigma=\chi_1\chi_2\cdots\chi_\tau$.

To distinguish between these two types of cycles when confusion is possible, we call the former kind an 
{\it edge cycle} and the latter kind a {\it permutation cycle}.  Notationally, we use Roman letters when 
discussing edge cycles and Greek letters when discussing permutation cycles.

%\medskip
\subsection{Permutation Expansion of the Determinant}

We recall that the {\it permutation expansion} of the determinant of an $n\times n$ matrix $M=(m_{ij})$ is the 
expansion of the determinant as
\begin{equation}
\det M = \sum_{\sigma \in S_{n}} (\sgn\, \sigma)m_{1,\sigma(1)}\cdots m_{n,\sigma(n)}.
\end{equation}
We will be considering non-zero terms in the permutation expansion of the determinant of the hamburger matrix 
$M_H$.  Because of the special block form of the hamburger matrix in Equation (\ref{BlockMatrix}), the 
permutations $\sigma$ that make non-zero contributions to this sum are products of disjoint cycles of either of 
two forms---the simple transposition
\begin{equation*}
\chi= (\phi_{11} ~ \omega_{11})
\label{TransCycle}
\end{equation*}
or the general permutation cycle
\begin{equation}
\chi= (\phi_{11} ~ \phi_{12} ~ \cdots ~ \phi_{1\mu_1} ~ \omega_{11} ~ \omega_{12} ~ \cdots ~ \omega_{1\nu_1} ~ 
\phi_{21} ~ \cdots\cdots ~ \phi_{\lam\mu_\lam} ~ \omega_{\lam 1} ~ \cdots ~ \omega_{\lam\nu_\lam}).
\label{PermCycle}
\end{equation}
In the first case, $\omega_{11}=\phi_{11}+k$.  In the second case, $1\leq \phi_{\i\kappa}\leq k$, $k+1\leq 
\omega_{\i\kappa}\leq 2k$, $\phi_{\i\kappa}<\phi_{\i,\kappa+1}$, and $\omega_{\i\kappa}>\omega_{\i,\kappa+1}$ for 
all $1\leq \i\leq \lam$ and relevant $\kappa$. The block matrix form also implies that 
$\phi_{\i\mu_\i}+k=\omega_{\i1}$, $\omega_{\i\nu_\i}-k=\phi_{\i+1,1}$ and $\omega_{\lam\nu_\lam}-k=\phi_{11}$.  
So that this permutation cycle is in standard form, we make sure that $\phi_{11}=\min_{\i,\kappa} 
\phi_{\i\kappa}$.  
In order to refer to this value later, we define a function $\Phi$ by $\Phi(\chi)=\phi_{11}$.
Each value $1\leq \phi_\i \leq k$ or $k+1\leq \omega_\i \leq 2k$ appears 
at most once for any $\sigma\in S_{2k}$.

We call a permutation cycle $\chi$ {\em minimal} if it is a transposition or if $\mu_\i=\nu_\i=2$ for all $\i$.  
This condition implies that we can write our general permutation cycles $\chi$ in the form
\begin{equation}
\chi= (\phi_{11} ~ \phi_{12} ~ \omega_{11} ~ \omega_{12} ~ \phi_{21} ~ \cdots ~ \phi_{\lam 2} ~ \omega_{\lam 1} ~ 
\omega_{\lam 2}),
\label{Minimal}
\end{equation}
with the same conditions as before.  We call a permutation $\sigma=\chi_1\cdots\chi_\tau$ {\em minimal} if each of 
its 
cycles $\chi_\iota$ is minimal.

%\medskip
\subsection{Walks Associated to a Permutation}

To each permutation $\chi\in S_{2k}$, we can associate one or more walks $c_\chi$ in $H$. 

\begin{figure}
\begin{center}
\begin{minipage}[c]{1in}
  $\chi = (1 ~ 2 ~ 3 ~ 6 ~ 4)$
\end{minipage}
\begin{minipage}[c]{2in}
  \epsfig{figure=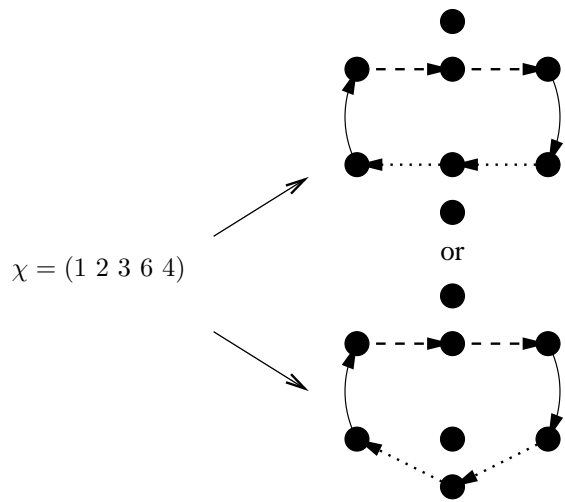}\\
\end{minipage}
\end{center}
  \caption{A permutation cycle $\chi$ and the two walks in $H$ associated to $\chi$}
  \label{CycleSystem}
\end{figure}

If $\chi$ is the transposition $\chi=(\phi_{11} ~ \omega_{11}),$ then we associate the 2-cycle 
$c_\chi:v_{\phi_{11}} \row w_{\omega_{11}} \row v_{\phi_{11}}$ to $\chi$.  To any permutation $\chi$ that is 
not a transposition, we can associate multiple walks $c_\chi$ by gluing together paths that follow $\chi$ in the 
following way.  If $\chi$ has the form 
of Equation (\ref{PermCycle}), then for each $1\leq i\leq \lam$, let $P_i$ be any path in $G_1$ that visits each 
of the vertices $v_{\phi_{i1}}$, $v_{\phi_{i2}}$, all the way through $v_{\phi_{i\mu_i}}$ in order.  Similarly, 
let $Q_i$ be any path in $G_2$ that visits each of the vertices $w_{\omega_{i1}}$, $w_{\omega_{i2}}$, through 
$w_{\omega_{i\nu_i}}$ in order.  For each choice of paths $P_i$ and $Q_i$, we have an additional possibility for 
the walk $c_\chi$; we can set
\begin{equation}
c_\chi: P_1 \stackrel{e_{\phi_{12}}}{\lrow} 
Q_1 \stackrel{e_{\phi_{21}}'}{\lrow} 
P_2 \stackrel{e_{\phi_{22}}}{\lrow} \cdots \cdots 
\stackrel{e_{\phi_{\lam 1}}'}{\lrow} 
P_{\lam} \stackrel{e_{\phi_{\lam 2}}}{\lrow} Q_{\lam}.
\label{Path}
\end{equation}
See Figure \ref{CycleSystem} for the choices of $c_{(12364)}$ in the hamburger graph presented in Figure 
\ref{NiceGraph}. We call $\lam$ the {\it number of $P$-paths} in $c_\chi$. The function $\Phi$, defined in the 
previous section, defines a partial ordering on walks in a walk system---we say that the associated walk 
$c_\chi$ comes before the associated walk $c_{\chi'}$ if $\Phi(\chi)<\Phi(\chi')$.  We call this the {\it 
initial term order}.  As in Section \ref{Sec:WHT}, we define the weight of a walk $c_\chi$ to be the product of 
the weights of all edges traversed by $c_\chi$.

%\medskip
\subsection{Walk System-Permutation Pairs}
\label{WSPPair}

We defined walk systems in Section 2, but we will see that the proof of Theorem \ref{HT} requires us to think of 
walk systems first as a permutation and second as a collection of walks determined by the permutation.  We will 
see that for cycle systems as presented initially, signs and weights are not changed by this recharacterization.

If $H$ is a hamburger graph with $k$ pairs of distinguished vertices, we define a walk system--permutation pair as 
follows.

{\definition A {\em walk system--permutation pair} (or {\em WSP-pair} for short) is a pair $(\C,\sigma)$, where 
$\sigma\in S_{2k}$ is a permutation and $\C$ is a collection of walks $c\in \C$ with the following property:  if 
the disjoint cycle representation of $\sigma$ is $\sigma=\chi_1\cdots\chi_\tau$, then $\C$ is a collection of 
$\tau$ 
walks $c_{\chi_\i}$, for $1\leq \i \leq \tau$, such that $c_{\chi_\i}$ is a walk associated to the 
permutation cycle $\chi_\i$.}

\medskip
\noindent 
We define the {\em weight} of a WSP-pair $(\C,\sigma)$ to be the product of the weights of the associated walks 
$c_\chi\in \C$.

The definition of WSP-pair implies that each permutation $\sigma$ yields many collections of walks $\C$, that 
collections of walks $\C$ may be associated to many permutations $\sigma$, but that any simple 
walk system $\C$ corresponds to one and only one minimal permutation $\sigma_m$.  This is because, given 
any path as in Equation (\ref{Path}), we can read off the initial and terminal vertices of each $P_i$ and $Q_i$ in 
order, producing a well-defined permutation cycle $\sigma_m$.  We define a WSP-pair $(\C,\sigma)$ to be {\it 
minimal} if $\sigma$ is a minimal permutation.

For a WSP-pair $(\C,\sigma)$, where $\sigma=\chi_1\cdots\chi_\tau$, we define the {\em sign} of the WSP-pair, 
$\sgn(\C,\sigma)$, to be $(-1)^{\ell}\sgn(\sigma),$ where $\lam_\chi$ is the number of $P$-paths in $c_\chi$ and 
where $\ell=\sum_{c_\chi\in \C} \lam_\chi$.  Alternatively, we could consider the sign of $(\C,\sigma)$ to be the 
product of the signs of its associated walks $c_\chi$, where the {\em sign} of $c_\chi$ is 
$\sgn(c_\chi)=(-1)^{\lam_\chi}\sgn(\chi)$. We say that a WSP-pair $(\C,\sigma)$ is {\it positive} if 
$\sgn(\C,\sigma)=+1$ and is {\it negative} if $\sgn(\C,\sigma)=-1$.

Note that if $(\C,\sigma)$ is a minimal WSP-pair, then $\sgn(c_\chi)=+1$ for a transposition $\chi$ and 
$\sgn(c_\chi)=(-1)^{\lam+1}$ if $\chi$ is of the form in Equation (\ref{Minimal}).  In particular, when 
$(\C,\sigma)$ is minimal and simple, its sign and weight is consistent with the definition given in the 
introduction.

%%%%%%%%%%%%%%%%%%%%%%%%%%%%%%%%%%%%%%%%%%%%%%%%%%%%%%%%%%%%%%%%
\section{Proof of the Hamburger Theorem}
\label{HamProof}

As mentioned in Section \ref{Outline}, we prove Lemmas \ref{Lem1} and \ref{Lem2} for the weighted version of 
the Hamburger Theorem, thereby proving Theorem \ref{WHT}. Theorem \ref{HT} follows as a special case of Theorem 
\ref{WHT}.

%%%%%%%%%%%%%%%%%%%%%%%%%%%%%%%%%%%%%%%%%%%%%%%%%%%%%%%%%%%%%%%%
\subsection{Proof of Lemma \ref{Lem1}, Part I}
\label{PfLem1}

Recall Properties 1 and 2 as well as Lemma \ref{Lem1}.

\medskip
\noindent
{\bf Property 1.} The walk system contains a walk that has a self-intersection.

\medskip
\noindent
{\bf Property 2.} The walk system has two intersecting walks, neither of which is a 2-cycle.

\medskip
\noindent
{\bf Lemma \ref{Lem1}.} {\em The set of all walk systems $\C$ that satisfy either Property 1 or Property 2 can be 
partitioned into equivalence classes, each of which contributes a net zero to the permutation expansion of the 
determinant of $M_H$.}

\begin{figure}
\begin{center}
  \epsfig{figure=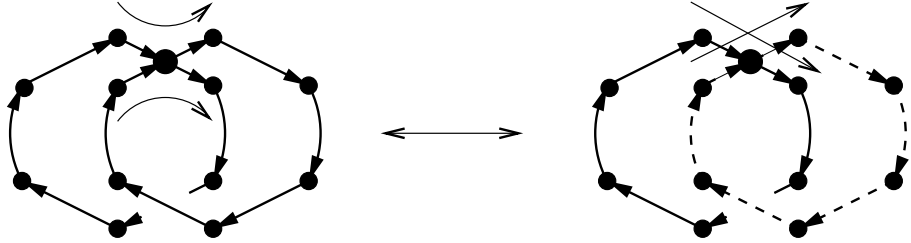}
  \caption{A self-intersecting cycle and its corresponding pair of intersecting cycles}
  \label{CycleSplitting}
\end{center}
\end{figure}

\medskip
The proof of Lemma \ref{Lem1} is a generalization of the involution principle, the idea of which comes from the 
picture presented in Figure \ref{CycleSplitting}.  Given a 
self-intersecting walk, changing the order of edge traversal at the vertex of self-intersection 
leads to breaking the one self-intersecting walk into two walks that intersect at that same vertex.  Since the 
edge set of the collection of walks has not changed, the weight of the two WSP-pairs is the same.  We have 
introduced a transposition into the sign of the permutation cycle of the WSP-pair; this changes the sign of the 
WSP-pair, so these two WSP-pairs will cancel in the permutation expansion of the determinant of $M_H$.

One can imagine that this means that to every self-intersecting WSP-pair we can associate one WSP-pair with two 
cycles intersecting.  However, more than one self-intersection may occur at this same point, and there may be 
additional walks that pass through that same point.  Exactly what this WSP-pair would cancel with is not clear.  
If we decide to break all the self-intersections so that we have some number $N$ of walks through our vertex, it 
is not clear how we should sew the cycles back together.  One starts to get the idea that we must consider all 
possible ways of sewing back together. Once we do just that, we have a family of WSP-pairs, all of the same 
weight, whose net contribution to the permutation expansion of the determinant is zero.

This idea is conceptually simple but the proof is notationally complicated.

\medskip
If $\C$ satisfies either Property 1 or Property 2, then there is some vertex of intersection, be it either a 
self-intersection or an intersection of two walks.  Our aim is to choose a well-defined {\em first} point 
of intersection at which we will build the family $\F$.  The initial term order gives an order on walks 
associated to permutation cycles; we choose the earliest walk $c_{\chi_\alpha}$ that has some vertex of 
intersection.  Once we have determined the earliest walk, we start at $v_{\Phi(c_{\chi_\alpha})}$ and follow 
the walk
\begin{equation*}
c_{\chi_a}: P_{1} \row Q_{1} \row P_{2} \row \cdots \row P_{m} \row Q_{m},
\end{equation*}
until we reach a vertex of intersection.  

In our discussion, we make the assumption that this first vertex of intersection is a vertex $v^*$ in $G_1$.  A 
similar argument exists if the first appearance occurs in $G_2$.  Notice that at $v^*$ there may be multiple 
self-intersections or multiple intersections of walks.  We will create a family $\F$ of WSP-pairs that takes into 
account each of these possibilities.

If we want to rigorously define the breaking of a self-intersecting walk at a vertex of self-intersection, 
we need to specify many different components of the WSP-pair $(C,\sigma)$.  First, we need to specify on which 
walk in $\C$ we are acting.  Next, we need to specify the vertex of self-intersection.  Since this 
self-intersection vertex may occur in multiple paths, we need to specify which two paths we interchange in the 
breaking process.

%%%%%%%%%%%%%%%%%%%%%%%%%%%%%%%%%%%%%%%%%%%%%%%%%%%%%%%%%%%%%%%%
\subsection{Definitions of Breaking and Sewing}

In the following paragraphs, we define ``breaking'' on WSP-pairs, which takes in a WSP-pair 
$(\C,\sigma)$, one of $\sigma$'s permutation cycles $\chi_\alpha$, the associated walk 
$c_{\chi_\alpha}$, paths $P_y$ and $P_z$ in $c_{\chi_\alpha}$, and the vertex $v^*$ in both $P_y$ and $P_z$ where 
$c_{\chi_\alpha}$ has a self-intersection. For simplicity, we assume that $v^*$ is not a distinguished vertex, but 
the argument still holds in that case.  The inverse of this operation is ``sewing''.

In this framework, 
$c_{\chi_\alpha}$ has the form
\begin{equation*}
c_{\chi_\alpha}:
P_1 \row Q_1 \row P_2 \row \cdots \row
P_y \row Q_y \row \cdots \row
P_z \row Q_z \row \cdots \row
P_\ell \row Q_\ell,
\end{equation*}
where the paths $P_y$ and $P_z$ in $G_1$ are separated into two halves as
\begin{equation*}
P_y : P_y^{(1)} \row v^* \row P_y^{(2)}
\end{equation*}
and
\begin{equation*}
P_z : P_z^{(1)} \row v^* \row P_z^{(2)}.
\end{equation*}

Remember that $P_y$ and $P_z$ are paths that stop over at various vertices depending on the permutation 
$\chi_\alpha$.  The vertex $v^*$ must have adjacent stop-over vertices in each of the two paths $P_y$ and $P_z$.  
Let the adjacent stop-over vertices in $P_y$ be $v_{\phi_y}$ and $v_{\phi_{y+1}}$ and the adjacent stop-over 
vertices in $P_z$ be $v_{\phi_z}$ and $v_{\phi_{z+1}}$.  

This implies $\chi_\alpha$ has the form
\begin{equation*}
\chi_\alpha = (
\phi_{11} ~ \cdots ~ \phi_{1\mu_1} ~ \omega_{11} ~ \cdots ~
\phi_{y} ~ \phi_{y+1} ~ \cdots ~
\phi_{z} ~ \phi_{z+1} ~ \cdots ~
\phi_{\lam\mu_\lam} ~ \omega_{\lam 1} ~ \cdots ~ \omega_{\lam\nu_\lam}).
\end{equation*}

We can now precisely define the result of breaking.  We define $\chi_\beta$ and $\chi_\gamma$ by 
splitting 
$\chi_\alpha$ as follows:
\begin{equation*}
\chi_\beta = (
\phi_{11} ~ \cdots ~ 
\phi_{y} ~ 
\phi_{z+1} ~ \cdots ~ \omega_{\lam\nu_\lam})
\end{equation*}
and
\begin{equation*}
\chi_\gamma = (\phi_{y+1} ~ \cdots ~ \phi_{z}),
\end{equation*}
with the necessary rewriting of $\chi_\gamma$ to have as its initial entry the value $\Phi(\chi_\gamma)$.  
Define $c_{\chi_\beta}$ and $c_{\chi_\gamma}$ to be 
\begin{equation*}
c_{\chi_\beta}:
P_1 \row Q_1 \row P_2 \row \cdots \row
P_y^{(1)} \row v^* \row 
P_z^{(2)}
 \row Q_z \row \cdots \row
P_\ell \row Q_\ell,
\end{equation*}
and
\begin{equation*}
c_{\chi_\gamma}:  
P_z^{(1)} \row v^* \row P_y^{(2)} \row Q_y \row \cdots \row Q_{z-1},
\end{equation*}
again changing the starting vertex of $c_{\chi_\gamma}$ to $v_{\Phi(c_{\chi_\gamma})}$.

We define the breaking of the WSP-pair with the above inputs to be the WSP-pair $(\C',\sigma')$ such that
\begin{equation*}
\C'=\C\cup\{c_{\chi_\beta},c_{\chi_\gamma}\}\setminus\{c_{\chi_\alpha}\} 
\end{equation*}
and 
\begin{equation*}
\sigma'=\sigma\chi_\alpha^{-1}\chi_\beta\chi_\gamma=\sigma\cdot (\phi_{y+1} ~ \phi_{z+1}).
\end{equation*}

The edge set of $\C$ is equal to the edge set of $\C'$, so the weight of the modified cycle 
systems is the same as the original.  Since we changed $\sigma$ to $\sigma'$ by multiplying only by a 
transposition, the sign of the modified WSP-pair is opposite to that 
of the original.

%%%%%%%%%%%%%%%%%%%%%%%%%%%%%%%%%%%%%%%%%%%%%%%%%%%%%%%%%%%%%%%%
%\medskip
\subsection{Proof of Lemma \ref{Lem1}, Part II}

Having defined breaking and sewing, we can continue the proof.

For any WSP-pair $(\C,\sigma)$ satisfying either Property 1 or Property 2, let $c_{\chi_\alpha}\in \C$ be the 
first 
walk in the initial term order with an vertex of intersection.  Let $v^*$ be the first vertex of 
intersection in $c_{\chi_\alpha}$.  Then for all walks $c$ with one or more self-intersections at $v^*$, 
continue to break $c$ at $v^*$ until there are no more self-intersections.  Define the resulting WSP-pair 
$(\C_s,\sigma_s)$ to be the {\it simple WSP-pair} associated to $(\C,\sigma)$.  In $(\C_s,\sigma_s)$, there is 
some number $N$ of general walks intersecting at vertex $v^*$.  There may be a 2-cycle intersecting $v^*$ as 
well, but this does not matter.

For any permutation $\xi\in S_N$, let $\xi=\zeta_1\zeta_2 \cdots \zeta_\eta$ be its cycle representation, where 
each $\zeta_\i$ is a cycle.  For each $1\leq \i \leq \eta$, sew together walks in order: if $\zeta_\i = 
(\delta_{\i 1} ~ \cdots ~ \delta_{\i \ep_\i})$, sew together $c_{\chi_{\delta_{\i 1}}}$ and $c_{\chi_{\delta_{\i 
2}}}$ at $v^*$.  Sew this result together with $c_{\chi_{\delta_{\i 3}}}$, and so on through $c_{\chi_{\delta_{\i 
\ep_\i}}}$.  Note that the result of these sewings is unique, and that every WSP-pair $(\C,\sigma)$ with 
$(\C_s,\sigma_s)$ as its simple WSP-pair can be obtained in this way, and no other WSP-pair appears.  We can 
perform this procedure for any $\xi\in S_N$; the sign of the resulting walk system $(C_\xi,\sigma_\xi)$ is 
$\sgn(C_\xi,\sigma_\xi)=\sgn(\xi)\,\sgn(\C_s,\sigma_s)$.  This means that the contribution to the determinant of 
the weights of all WSP-pairs in the family $\F$ is
\begin{equation}
\sum_{\xi\in S_N} \sgn(\xi)\,\sgn(\C_s,\sigma_s)\,\wt(\C_s,\sigma_s)
=\sgn(\C_s,\sigma_s)\,\wt(\C_s,\sigma_s) \sum_{\xi\in S_N} \sgn(\xi) = 0.
\end{equation}
So elements of the same family cancel out in the determinant of $M_H$, giving that the contributions of all 
WSP-pairs satisfying either Property 1 or Property 2 cancel out in the hamburger determinant. \hfill$\bullet$

%%%%%%%%%%%%%%%%%%%%%%%%%%%%%%%%%%%%%%%%%%%%%%%%%%%%%%%%%%%%%%%%
\subsection{Proof of Lemma \ref{Lem2}}
\label{PfLem2}

Recall Properties 3 and 4 as well as Lemma \ref{Lem2}.

\medskip
\noindent
{\bf Property 3.} The walk system has two intersecting walks, one of which is a 2-cycle.

\medskip
\noindent
{\bf Property 4.} The walk system is not minimal.

\medskip
\noindent
{\bf Lemma \ref{Lem2}.} {\em The set of all walk systems $\C$ that satisfy neither Property 1 nor Property 2 and 
that satisfy either Property 3 or Property 4 can be partitioned into equivalence classes, each of which 
contributes a net zero to the permutation expansion of the determinant of $M_H$.}

\medskip
In the proof of Lemma \ref{Lem2}, we do not base our family $\F$ around a singular vertex; instead, we find a set 
of violations that each member of the family has.  If the WSP-pair satisfies the hypotheses of Lemma \ref{Lem2}, 
then either there is some 2-cycle $c:v_\phi\row w_\omega\row v_\phi$ that intersects with some other walk or 
there is some non-minimal permutation cycle.

Define a set of indices $I\subseteq [k]$ of violations, of which an integer can become a member in one of two 
ways.  If $(\C,\sigma)$ is not minimal, there is at least one permutation cycle $\chi_\alpha$ with more than two 
consecutive $\phi$'s or $\omega$'s in its cycle notation.  For any intermediary $\iota$ between two $\phi$'s or 
$\iota+k$ between two $\omega$'s, place $\iota$ in $I$.  For example, if $\chi_\alpha=(\cdots ~ \phi' ~ \i ~ 
\phi'' ~ \cdots)$, we place $\i\in I$.  Alternatively, there may be a 2-cycle $c: v_i \row w_{k+i} \row v_i$ such 
that either $v_i$ is a vertex in some other walk $c_{\chi_\beta}$ or $w_{k+i}$ is a vertex in some other cycle 
$c_{\chi_\gamma}$, or both.  We declare this $i$ to be in $I$ as well.

Note that any WSP-pair $(\C,\sigma)$ satisfying either Property 3 or Property 4 has a non-empty set $I$.  From our 
original WSP-pair, associate a minimal WSP-pair $(\C_m,\sigma_m)$ by removing any transposition $\chi_\alpha$ from 
$\sigma$ and its corresponding 2-cycle $c_{\chi_\alpha}$ from $\C$, and also for any non-minimal permutation cycle 
$\chi_\beta$ we remove any intermediary $\phi$'s or $\omega$'s from $\sigma$.  We do not change the associated 
walk $c_{\chi_\beta}$ in $\C$ since it still corresponds to this minimized permutation cycle.

Let $i$ be any element in $I$.  Since $i\in I$, the 2-cycle $c_i: v_i \row w_{k+i} \row v_i$ intersects some walk 
of $\C$ either at $v_i$, at $w_{k+i}$, or both.  So there are four cases:

\medskip
\noindent {\bf Case 1.} $c_i$ intersects a walk $c_{\chi_\beta}$ at $v_i$ and no walk at $w_{k+i}$.

\medskip
\noindent {\bf Case 2.} $c_i$ intersects a walk $c_{\chi_\gamma}$ at $w_{k+i}$ and no walk at $v_i$.

\medskip
\noindent {\bf Case 3.} $c_i$ intersects a walk $c_{\chi_\beta}$ at $v_i$ and the same walk again at 
$w_{k+i}$.

\medskip
\noindent {\bf Case 4.} $c_i$ intersects a walk $c_{\chi_\beta}$ at $v_i$ and another walk 
$c_{\chi_\gamma}$ at $w_{k+i}$.

\begin{figure}
\begin{center}
  \epsfig{figure=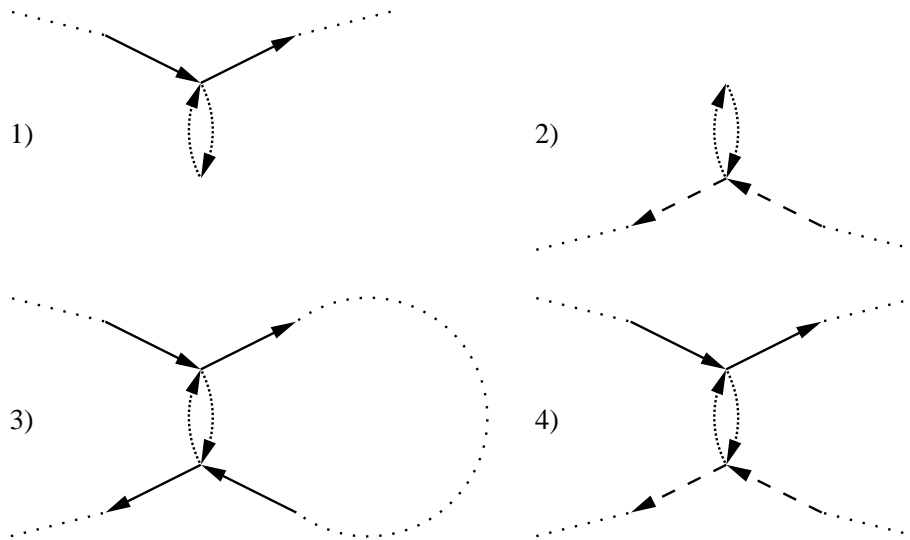}
  \caption{The four cases in which a walk can intersect with a 2-cycle}
  \label{2Cycles}
\end{center}
\end{figure}

\medskip
\noindent
See Figure \ref{2Cycles} for a picture.

In order to build the family $\F$, we modify the minimal WSP-pair $(\C_m,\sigma_m)$ at each $i$ with some choice 
of options, each giving a possible WSP-pair that has $(\C_m,\sigma_m)$ as its minimal WSP-pair.  We define a set 
of two or four options $\O_i$ for each $i$.

In Case 1, there are two options.  Let $o_1$ be the option of including the walk $c_i$ in $\C_m$ and its 
corresponding transposition $\chi_\alpha$ in $\sigma_m$.  Let $o_2$ be the option of including the intermediary 
$\phi=i$ in $\chi_\beta$ in the position where $c_{\chi_\beta}$ passes through $v_i$.

[Note that we could not apply both options at the same time since the resulting $\chi_\alpha$ and $\chi_\beta$ 
would not be disjoint and therefore is not a term in the permutation expansion of the determinant.]

Similarly in Case 2, there are two options.  Let $o_1$ be the option of including the walk $c_i$ in $\C_m$ and its 
corresponding transposition $\chi_\alpha$ in $\sigma_m$.  Let $o_2$ be the option of including the intermediary 
$\omega=i+k$ in $\chi_\gamma$ in the position where $c_{\chi_\gamma}$ passes through $w_{k+i}$.

In Case 3, there are four options.  Let $o_1$ be the option of including the walk $c_i$ in $\C_m$ and its 
corresponding transposition $\chi_\alpha$ in $\sigma_m$.  Let $o_2$ be the option of including the intermediary 
$\phi=i$ in $\chi_\beta$ in the position where $c_{\chi_\beta}$ passes through $v_i$.  Let $o_3$ be the option of 
including the intermediary $\omega=i+k$ in $\chi_\beta$ in the position where $c_{\chi_\beta}$ passes through 
$w_{k+i}$.  Let $o_4$ be the option of including both intermediaries $\phi=i$ and $\omega=i+k$ in $\chi_\beta$ in 
the respective positions where $c_{\chi_\beta}$ passes through $v_i$ and $w_{k+i}$.

In Case 4, there are four options.  Let $o_1$ be the option of including the walk $c_i$ in $\C_m$ and its 
corresponding transposition $\chi_\alpha$ in $\sigma_m$.  Let $o_2$ be the option of including the intermediary 
$\phi=i$ in $\chi_\beta$ in the position where $c_{\chi_\beta}$ passes through $v_i$.  Let $o_3$ be the option of 
including the intermediary $\omega=i+k$ in $\chi_\gamma$ in the position where $c_{\chi_\gamma}$ passes through 
$w_{k+i}$.  Let $o_4$ be the option of including intermediary $\phi=i$ in $\chi_\beta$ in the position where 
$c_{\chi_\beta}$ passes through $v_i$ and intermediary $\omega=i+k$ in $\chi_\gamma$ in the position where 
$c_{\chi_\gamma}$ passes through $w_{k+i}$.

\begin{figure}
\begin{center}
  \epsfig{figure=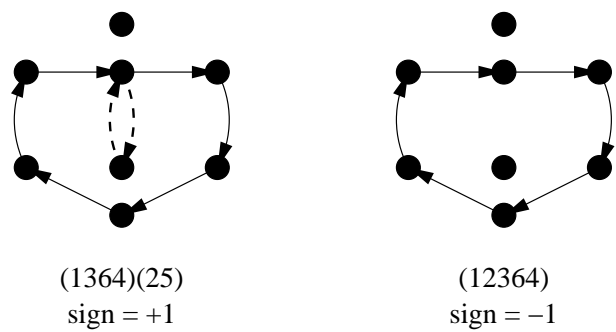}
  \caption{A family of two walk systems that cancel in the hamburger determinant}
  \label{Lem2Example}
\end{center}
\end{figure}

\medskip 
Corresponding to the associated minimal WSP-pair $(\C_m,\sigma_m)$ and index set $I$, define the family $\F$ to be 
the set of WSP-pairs $(\C_f,\sigma_f)$ by exercising all combinations of options $o_{j_i}\in \O_i$ on 
$(\C_m,\sigma_m)$ for each $i\in I$.  Note that every WSP-pair derived in this fashion is a WSP-pair satisfying 
the hypotheses of Lemma \ref{Lem2} and is such that its associated minimal WSP-pair is $(\C_m,\sigma_m)$.  There 
is also no other WSP-pair $(\C',\sigma')$ with $(\C_m,\sigma_m)$ as its minimal WSP-pair.  See Figure 
\ref{Lem2Example} for a canceling family of two walk systems, one of which is the non-minimal 
walk system from Figure \ref{CycleSystem}.

Every WSP-pair in $\F$ has the same weight since each option changes the edge set of $\C$ by at most a 2-cycle 
$v_i\row w_{k+i} \row v_i$, where $2\leq i \leq k-1$ and each of those 2-cycles contributes a multiplier of $1$ to 
the weight of the WSP-pair.  Note that the peculiar bounds in this restriction are due to the fact that no 
permutation $\sigma$ can ever include two cycles through any of the vertices $v_1$, $w_{k+1}$, $v_k$, or $w_{2k}$.  
The sign of every WSP-pair in $\F$ is determined by the set of options applied to $(\C_m,\sigma_m)$.  Each option 
$o_j$ contributes a multiplicative $\pm 1$ depending on $j$---$\sgn(o_1)=+1$, $\sgn(o_2)=-1$, and $\sgn(o_3)=-1$ 
and $\sgn(o_4)=+1$ if they exist.  

Each family $\F$ contributes a cumulative weight of 0.  This is because there are the same number of positive 
$(\C_f,\sigma_f)\in\F$ as negative $(\C_f,\sigma_f)\in F$.  We see this since if $i^*\in I$ then any set of 
WSP-pairs that exercise the same options $o_{j_i}$ for all $i\ne i^*$ has either two or four members (depending on 
the case of $i^*$), which split evenly between positive and negative WSP-pairs.  Since each family contributes a 
net zero to the hamburger determinant, the total contribution from WSP-pairs satisfying the hypotheses of Lemma 
\ref{Lem2} is zero.\hfill$\bullet$

\medskip
Since every WSP-pair appears once in the permutation expansion of $\det M_H$, the only WSP-pairs that contribute 
to the sum are those that are simple and minimal.  Therefore $\det M_H$ is the sum over such WSP-pairs of 
$\sgn(\C,\sigma)\cdot\wt(\C,\sigma)$.  This establishes Theorem \ref{WHT} and gives Theorem \ref{HT} as a special 
case.

%%%%%%%%%%%%%%%%%%%%%%%%%%%%%%%%%%%%%%%%%%%%%%%%%%%%%%%%%%%%%%%%
\section{Applications of the Hamburger Theorem}
\label{Apps}

We discuss first the application of the Hamburger Theorem in the case where the region is an Aztec diamond, 
mirroring results of Brualdi and Kirkland.  Then we discuss the results from the case where the region is an Aztec 
pillow, and along the way, we explain how to implement the Hamburger Theorem when our region is any generalized 
Aztec pillow.

%%%%%%%%%%%%%%%%%%%%%%%%%%%%%%%%%%%%%%%%%%%%%%%%%%%%%%%%%%%%%%%%
\subsection{Domino Tilings and Digraphs}

\begin{figure}
\begin{center}
\begin{tabular}{ccc}
  \epsfig{figure=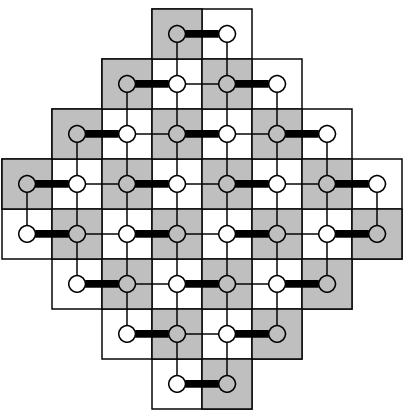} &
  \epsfig{figure=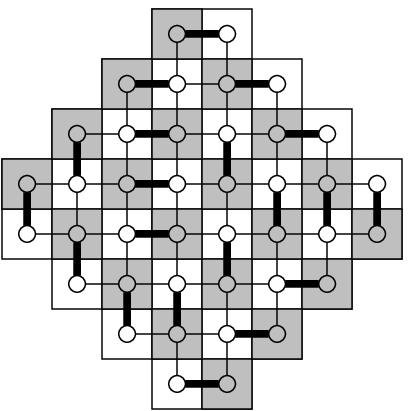} &
  \epsfig{figure=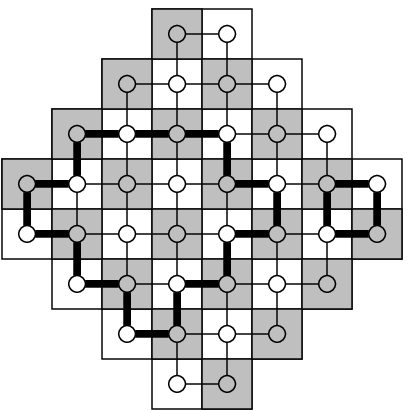} 
\end{tabular}
  \caption{The symmetric difference of two matchings gives a cycle system}
  \label{AD}
\end{center}
\end{figure}

The Hamburger Theorem applies to counting the number of domino tilings of Aztec diamonds and generalized Aztec 
pillows.  To illustrate this connection, we count domino tilings of the Aztec diamond by counting an equivalent 
quantity, the number of matchings on the dual graph $G$ of the Aztec diamond.  We use the natural matching $N$ of 
horizontal neighbors in $G$, as exemplified in Figure \ref{AD}a on $AD_4$, as a reference point.  Given any other 
matching $M$ on $G$, such as in Figure \ref{AD}b, its symmetric difference with the natural matching is a 
collection of cycles in the graph, such as in Figure \ref{AD}c.

When we orient the edges in $N$ from black vertices to white vertices and 
in $M$ from white vertices to black vertices, the symmetric difference 
becomes a collection of directed cycles.  Notice that edges in the upper half 
of $G$ all go from left to right and the edges in the bottom half of $G$ 
go from right to left.

\begin{figure}
\begin{center}
\begin{tabular}{cc}
  \epsfig{figure=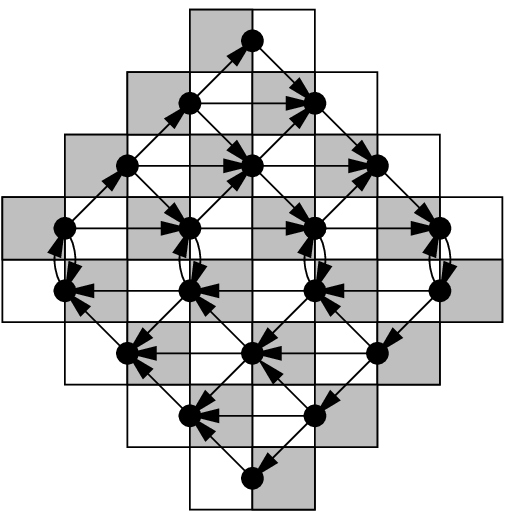} &
  \epsfig{figure=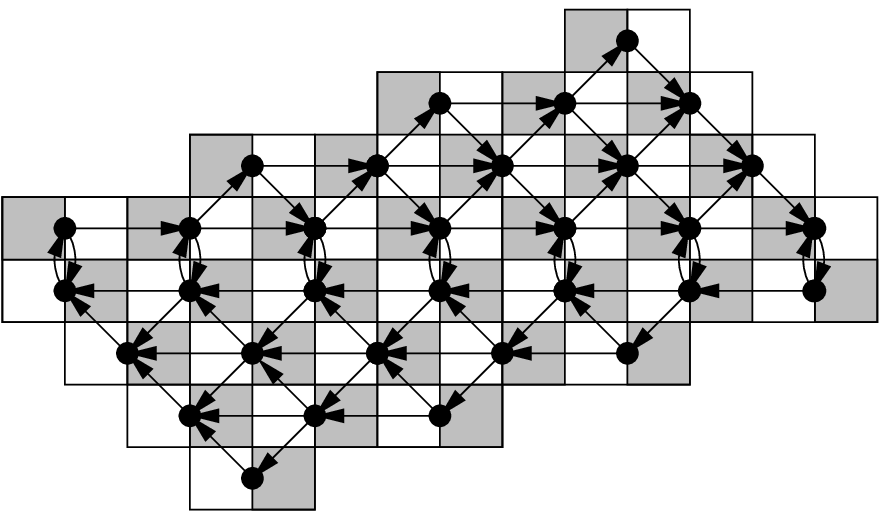}
\end{tabular}
  \caption{The hamburger graph for an Aztec diamond and an Aztec pillow}
  \label{ADAPHamGraph}
\end{center}
\end{figure}

Since the edges of $N$ always appear in the cycles, we can contract the 
edges of $N$ to points; then the graph retains its structure in terms of 
the cycle systems it produces. This new graph $H$ is of the form in Figure 
\ref{ADAPHamGraph}a.
This argument shows that the number of domino tilings of an Aztec diamond 
equals the number of cycle systems of this new condensed graph, called the 
region's {\em digraph}.  

We wish to concretize this notion of a {\em digraph} of the Aztec diamond $AD_n$. Given the natural tiling of an 
Aztec diamond consisting solely of horizontal dominoes, we place a vertex in the center of every domino.  The 
edges of this digraph fall into three families.  From every vertex in the top half of the diamond, create edges to 
the east, to the northeast, and to the southeast whenever there is a vertex there.  From every vertex in the 
bottom half of the diamond, form edges to the west, to the southwest, and to the northwest whenever there is a 
vertex there.  Additionally, label the bottom vertices in the top half $v_1$ through $v_n$ from west to east and 
the top vertices in the bottom half $w_{n+1}$ through $w_{2n}$, also west to east.  For all $i$ between $1$ and 
$n$, create a directed edge from $v_i$ to $w_{n+i}$ and from $w_{n+i}$ to $v_i$.  The result when this 
construction is applied to $AD_4$ is a graph of the form in Figure \ref{ADAPHamGraph}a.  By construction, this 
graph is a hamburger graph.

{\theorem The digraph of an Aztec diamond is a hamburger graph.}

\medskip
Since both the upper half of the digraph and the lower half of the digraph 
are strongly planar, there are no negative cycle systems.  This implies that 
the determinant of the corresponding hamburger matrix counts exactly the 
number of cycle systems in the digraph.

{\corollary The number of domino tilings of an Aztec diamond is the determinant of its hamburger matrix.}

%%%%%%%%%%%%%%%%%%%%%%%%%%%%%%%%%%%%%%%%%%%%%%%%%%%%%%%%%%%%%%%%
\subsection{Applying the Hamburger Theorem to Aztec Diamonds}
\label{applyAD}

\begin{figure}
\begin{center}
  \epsfig{figure=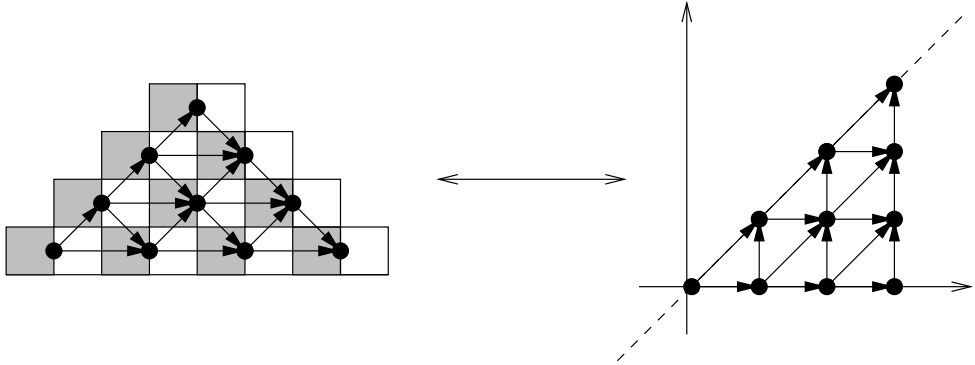}
  \caption{The equivalence between paths in $D$ and lattice paths in the first quadrant}
  \label{LatticePaths}
\end{center}
\end{figure}

In order to apply Theorem \ref{HT} to count the number of tilings of $AD_n$, we need to find the number of paths 
in the upper half of $D$ from $v_i$ to $v_j$ and the number of paths in the lower half of $D$ from $w_{n+j}$ to 
$w_{n+i}$.  The key observation is that by the equivalence in Figure \ref{LatticePaths}, we are in effect counting 
the number of paths from $(i,i)$ to $(j,j)$ using steps of size $(0,1)$, $(1,0)$, or $(1,1)$ that do not pass 
above the line $y=x$.  This is a combinatorial interpretation for the $(j-i)$-th large Schr\"oder number.  The 
large Schr\"oder numbers $s_0,\hdots,s_5$ are $1, 2, 6, 22, 90, 394$.  This is sequence number A006318 in the 
Encyclopedia of Integer Sequences \cite{EIS}.

{\corollary The number of domino tilings of the Aztec 
diamond $AD_n$ is equal to $$\#AD_n=\det \left[ \begin{array}{cc} S_n & I_n \\ 
-I_n & S_n^T \\ \end{array} \right],$$ where $S_n$ is an upper triangular 
matrix with the $i$-th large Schr\"oder number on its $i$th superdiagonal.}

\medskip  For example, when $n=6$, the matrix $S_6$ is
$$S_6 = \left[ 
\begin{array}{cccccc} 
1 & 2 & 6 & 22 & 90 & 394 \\
0 & 1 & 2 & 6 & 22 & 90 \\
0 & 0 & 1 & 2 & 6 & 22 \\
0 & 0 & 0 & 1 & 2 & 6 \\
0 & 0 & 0 & 0 & 1 & 2  \\
0 & 0 & 0 & 0 & 0 & 1   \\
\end{array} \right].$$ 

\medskip 
Brualdi and Kirkland prove a similar determinant formula for the number of tilings of an Aztec diamond in a 
matrix-theoretical fashion based on the $n(n+1)\times n(n+1)$ Kasteleyn matrix of the graph $H$ and a Schur 
complement calculation.  The Hamburger Theorem gives a purely combinatorial way to reduce the calculation of the 
number of tilings of an Aztec diamond to the calculation of a $2n\times 2n$ Hamburger determinant.  

Following cues from Brualdi and Kirkland, we can reduce this to an $n \times n$ determinant via a Schur 
complement.  For uses and the history of the Schur Complement, see the new book \cite{Zhang} by Zhang. In the case 
of the block matrix $M_H$ in Equation (\ref{BlockMatrix}), taking the Schur complement of $B$ in $M_H$ gives
\begin{equation}
\det M_H = \det B \cdot\det (A+D_1B^{-1}D_2)=\det(A+D_1B^{-1}D_2),
\end{equation}
since $B$ is a lower triangular matrix with 1's on the diagonal.  In this way, every hamburger determinant can be 
reduced to a smaller determinant of a Schur complement matrix.  In the case of a simple hamburger graph where 
$D_2=D_1=I$, the determinant reduces further to $\det(A+B^{-1})$.  Lastly, in the case where the 
hamburger graph is rotationally symmetric, $B=JAJ$ where $J$ is the exchange matrix, which consists of 1's down 
the skew-diagonal 
and 0's elsewhere.  This implies that we can write the determinant in terms of just the submatrix $A$, i.e., 
$\det(A+JA^{-1}J)$. (Note that $J^{-1}=J$.)

{\corollary The number of domino tilings of the Aztec diamond $AD_n$ is equal to $\det (S_n + J_nS_n^{-1}J_n)$, where $J_n$ is the 
$n\times n$  exchange 
matrix.}

\medskip
In the case of a hamburger graph $H$, we call this Schur complement, $A+B^{-1}$, a {\em reduced hamburger matrix}.
In the Aztec diamond graph example above, we can thus calculate the number of tilings of the Aztec diamond $AD_6$ as 
follows.  The inverse of $S_6$ is 
$$S_6^{-1} = \left[ 
\begin{array}{rrrrrr} 
1 & -2 & -2 & -6 & -22 & -90  \\
0 & 1 & -2 & -2 & -6 & -22 \\
0 & 0 & 1 & -2 & -2 & -6  \\
0 & 0 & 0 & 1 & -2 & -2 \\
0 & 0 & 0 & 0 & 1 & -2  \\
0 & 0 & 0 & 0 & 0 & 1   \\
\end{array} \right],$$ 
which implies that the determinant of the reduced hamburger matrix
$$M_6 = \left[ 
\begin{array}{rrrrrr} 
2 & 2 & 6 & 22 & 90 & 394 \\
-2 & 2 & 2 & 6 & 22 & 90 \\
-2 & -2 & 2 & 2 & 6 & 22 \\
-6 & -2 & -2 & 2 & 2 & 6 \\
-22 & -6 & -2 & -2 & 2 & 2  \\
-90 & -22 & -6 & -2 & -2 & 2   \\
\end{array} \right]$$ 
gives the number of tilings of $AD_6$.  

Brualdi and Kirkland were the first to find such a determinantal formula for the number of tilings of an Aztec 
Diamond \cite{BruKirk}.  Their matrix was different only in the fact that each entry was multiplied by $(-1)$ and 
there was a multiplicative factor of $(-1)^n$.  Brualdi and Kirkland were able to calculate the sequence of 
determinants, $\{\det M_n\}$, using a $J$-fraction expansion, which only works when matrices are Toeplitz or 
Hankel.

%%%%%%%%%%%%%%%%%%%%%%%%%%%%%%%%%%%%%%%%%%%%%%%%%%%%%%%%%%%%%%%%
\subsection{Applying the Hamburger Theorem to Aztec Pillows}
\label{applyAP}

\begin{figure}
\begin{center}
  \epsfig{figure=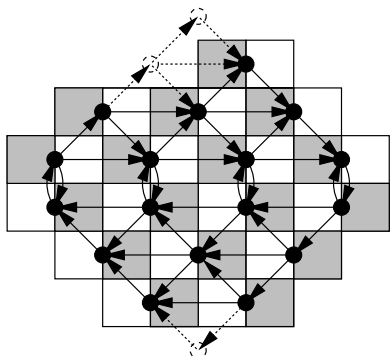}
  \caption{The digraph of a generalized Aztec pillow from the digraph of an Aztec diamond}
  \label{GenAPdigraph}
\end{center}
\end{figure}

\medskip 
A generalized Aztec pillow can be produced by restricting the placement of certain dominoes in an Aztec diamond in 
order to generate the desired boundary. We define the digraph of a generalized Aztec pillow to be the restriction 
of the digraph of an Aztec diamond to the vertices on the interior of the pillow.  An example is presented Figure 
\ref{GenAPdigraph}.  Since the generalized Aztec pillow's digraph is a restriction of the Aztec diamond's 
digraph, we have the following corollary.

{\corollary The digraph of an Aztec pillow or a generalized Aztec pillow is a hamburger graph.}

\medskip
Aztec pillows were introduced in part because of an intriguing conjecture for the number of their tilings given by 
Propp \cite{Propp}.

\medskip
\noindent
{\conjecture[Propp's Conjecture]  
The number of tilings of an Aztec pillow $AP_n$ is a larger number squared times a smaller number.  We write 
$\#AP_n=\ell_n^2s_n$. In addition, depending on the parity of $n$, the smaller number $s_n$ satisfies a simple 
generating function. For $AP_{2m}$, the generating function is
\begin{equation*}
\sum_{m=0}^\infty s_{2m}x^m = {(1+3x+x^2-x^3)}/{(1-2x-2x^2-2x^3+x^4)}.
\label{GF1}
\end{equation*}
For $AP_{2m+1}$, the generating function is
\begin{equation*}
\sum_{m=0}^\infty s_{2m+1}x^m = {(2+x+2x^2-x^3)}/{(1-2x-2x^2-2x^3+x^4)}.
\label{GF2}
\end{equation*}
}

\medskip
\noindent 
My hope for the Hamburger Theorem was that it would give a proof of Propp's Conjecture.  We 
shall see that although we achieve a faster determinantal method to calculate $\#AP_n$, the sequence of 
determinants is not calculable by known methods.

\begin{figure}
\begin{center}
  \epsfig{figure=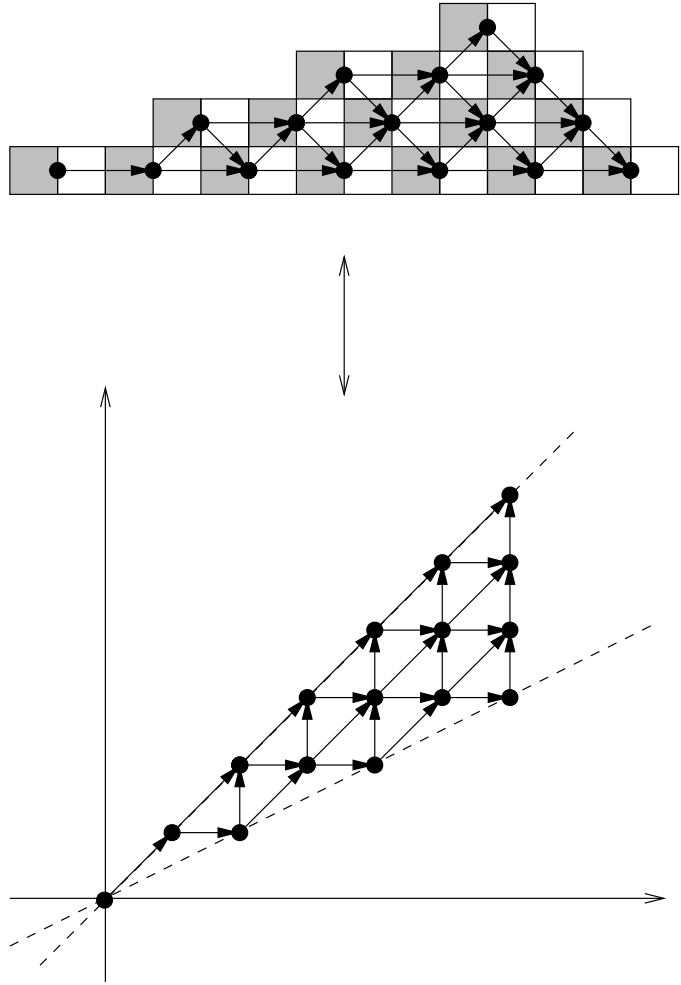}
  \caption{The equivalence between paths in $D$ and lattice paths in the first quadrant}
  \label{APLatticePaths}
\end{center}
\end{figure}

Using the same method as for Aztec diamonds, creating the hamburger graph $H$ for an Aztec pillow gives Figure 
\ref{ADAPHamGraph}b.  Counting the number of paths from $v_i$ to $v_j$ and from $w_{n+j}$ to $w_{n+i}$ in 
successively larger Aztec pillows gives us the infinite upper-triangular array $S=(s_{i,j})$ of {\em modified 
Schr\"oder numbers} defined by the following combinatorial interpretation.  Let $s_{i,j}$ be the number of paths 
from $(i,i)$ to $(j,j)$ using steps of size $(0,1)$, $(1,0)$, and $(1,1)$, not passing above the line $y=x$ nor 
below the line $y=x/2$.  The equivalence between these paths and paths in the upper half of the hamburger graph 
is shown in Figure \ref{APLatticePaths}.
The principal $7\times 7$ minor matrix $S_7$ of $S$ is 
$$S_7=\left[
\begin{array}{rrrrrrr}
1 & 1 & 2 & 5 & 16 & 57 & 224\\
0 & 1 & 2 & 5 & 16 & 57 & 224\\
0 & 0 & 1 & 2 &  6 & 21 & 82\\
0 & 0 & 0 & 1 &  2 &  6 & 22\\
0 & 0 & 0 & 0 &  1 &  2 & 6\\
0 & 0 & 0 & 0 &  0 &  1 & 2\\
0 & 0 & 0 & 0 & 0 &  0 &  1\\
\end{array}
\right].$$
This table of numbers is sequence number A114292 in the Encyclopedia of Integer Sequences \cite{EIS}.

{\theorem The number of domino tilings of an Aztec pillow of order $n$ is equal to 
$$\#AP_n=\det\left[ \begin{array}{cc} S_n & I_n \\ -I_n & J_nS_nJ_n \\ \end{array} \right],$$ where $S_n$ is the $n\times n$ 
principal 
submatrix of $S$ and $J_n$ is the $n\times n$ exchange matrix.}

\medskip
As in the case of Aztec diamonds, we can calculate the 
reduced hamburger matrix through a Schur calculation.  The inverse of $S_6$ is 
$$S_6^{-1}=\left[
\begin{array}{rrrrrrr}
1 & -1 & 0 & 0 &  0 &  0 \\
0 & 1 & -2 & -1 &  -2 &  -5 \\
0 & 0 & 1 & -2 &  -2 &  -5 \\
0 & 0 & 0 & 1 &  -2 &  -2 \\
0 & 0 & 0 & 0 &  1 &  -2 \\
0 & 0 & 0 & 0 &  0 &  1 \\
\end{array}
\right]$$
and the resulting reduced hamburger matrix for $AP_6$ is
$$M_6=\left[
\begin{array}{rrrrrrr}
 2 &  1 &  2 & 5 & 16 & 57 \\
-2 &  2 &  2 & 5 & 16 & 57 \\
-2 & -2 &  2 & 2 &  6 & 21 \\
-5 & -2 & -2 & 2 &  2 &  6 \\
-5 & -2 & -1 &-2 &  2 &  2 \\
 0 &  0 &  0 & 0 & -1 &  2 \\
\end{array}
\right].$$

This gives us a much faster way to calculate the number of domino tilings of an Aztec pillow than was known 
previously.  We have reduced the calculation of the $O(n^2)\times O(n^2)$ Kasteleyn--Percus determinant to an 
$n\times n$ reduced hamburger matrix.  To be fair, the Kasteleyn--Percus matrix has $-1$, $0$, and $+1$ entries 
while the reduced hamburger matrix may have very large entries, which makes theoretical running time comparisons 
difficult.  Experimentally, when calculating the number of domino tilings of $AP_{14}$ using Maple 8.0 on a 447 
MHz Pentium III processor, the determinant of the $112\times 112$ Kasteleyn--Percus matrix took 25.3 seconds 
while the determinant of the $14\times 14$ reduced hamburger matrix took less than 0.1 seconds.

Whereas we have a very understandable determinantal formula for the number of tilings of the region, this does 
not translate into a proof of Propp's Conjecture because we cannot calculate the determinant of the 
matrices $M_n$ explicitly.  We cannot apply a $J$-fraction expansion as Brualdi and Kirkland did since the 
reduced hamburger matrix is not Toeplitz or Hankel.

\section{Generalizations}
\label{Generalizations}
In this section we present a counterexample to a natural generalization of the Hamburger Theorem and discuss an 
intriguing generalization of Propp's Conjecture.

%%%%%%%%%%%%%%%%%%%%%%%%%%%%%%%%%%%%%%%%%%%%%%%%%%%%%%%%%%%%%%%%
\subsection{Counterexample to Generalizations of the Hamburger Theorem}

The structure of the hamburger graph presented in Figure \ref{Hamburger} seems restrictive, so the question 
naturally arises whether it is somehow necessary.  Can the edge set $E_3$ between graphs $G_1$ and $G_2$ be an 
arbitrary bipartite graph?  The answer in general is no.  

\begin{figure}
\begin{center}
  \epsfig{figure=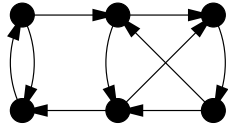}
  \caption{A simple generalized hamburger graph}
  \label{CtrEx1}
\end{center}
\end{figure}

\begin{figure}
\begin{center}
  \epsfig{figure=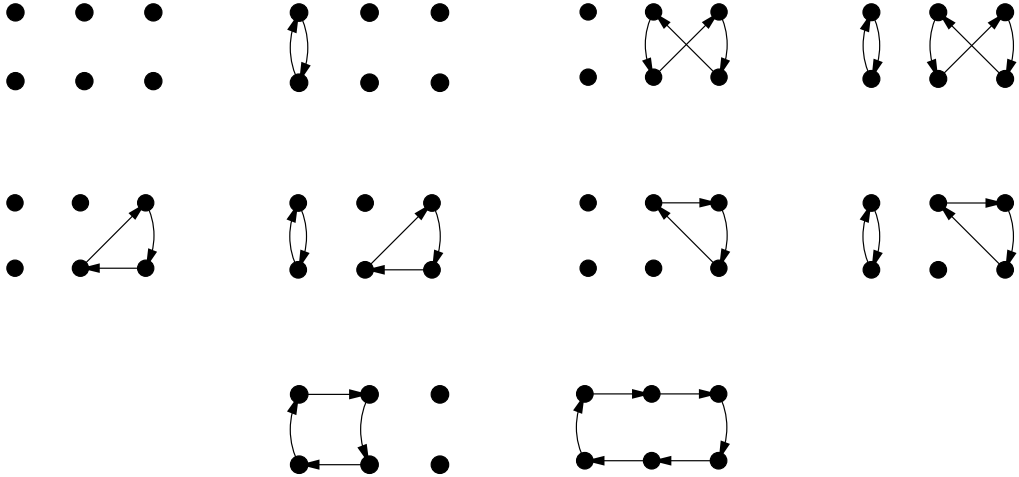}
  \caption{The ten cycle systems for the generalized hamburger graph in Figure \ref{CtrEx1}}
  \label{CtrEx2}
\end{center}
\end{figure}

Take for example the simple graph $H$ in Figure 
\ref{CtrEx1}.
As one can count by hand, there are 10 distinct cycle systems in 
$H$, as enumerated in Figure \ref{CtrEx2}.
Creating the hamburger matrix that would correspond to this graph gives
\begin{equation*}
M_H=
\left[
\begin{array}{cccccc}
 1 & 1 & 1 & 1 & 0 & 0 \\
 0 & 1 & 1 & 0 & 1 & 0 \\
 0 & 0 & 1 & 0 & 0 & 1 \\
-1 & 0 & 0 & 1 & 0 & 0 \\
0 & 0 & -1 & 1 & 1 & 0 \\
0 & -1 & 0 & 1 & 1 & 1 \\
\end{array}
\right],
\end{equation*}
where the lower left block describes paths from $w_4$, $w_5$, and $w_6$ to $v_1$, $v_2$, and $v_3$ and we have 
changed no other block of the matrix. This seems the most logical extension of the hamburger matrix.

The determinant of this matrix is $-5$.  Even though one might expect to need a new sign convention on cycle 
systems in generalized hamburger graphs, any sign convention would necessarily conserve the parity of the number 
of cycle systems.  Therefore there is no $+1$/$-1$ labeling of the cycle systems in Figure \ref{CtrEx2} that would 
allow $\det M_H=\sum_{c\in\CC} \sgn(c)$.  This means that either the matrix for this generalized hamburger graph 
is not correct or that the theorem does not hold in general.

The only hope for my idea of how to generalize the hamburger matrix is what happens when one considers the 
difference between even and odd permutations.  
In the counterexample above, the edge set $E_3$ sends vertex $w_4$ to vertex $v_1$, $w_5$ to $v_3$, and $w_6$ to 
$v_2$.  This can be thought of as an odd permutation of a hamburger's edge set (sending $w_4$ to $v_1$, $w_5$ to 
$v_2$, and $w_6$ to $v_3$).  In my limited calculations so far, whenever the permutation of the natural edge set 
is even, the parity of the number of cycle systems, counted by hand, is equal to the parity of the determinant of 
the matrix.  I have not found a sign convention for these graphs that explains the calculated determinant, but 
being equal in parity allows for hope of a possible extension to the Hamburger Theorem.

%%%%%%%%%%%%%%%%%%%%%%%%%%%%%%%%%%%%%%%%%%%%%%%%%%%%%%%%%%%%%%%%
\subsection{Generalizing Propp's Conjecture}

A key motivation in the study of Aztec pillows is Propp's Conjecture for 
3-pillows, presented in Section \ref{applyAP}.  Through experimental 
calculations, I have expanded and extended the conjecture.

For one, I extend the conjecture to odd pillows.  For any odd number $q$, we can use the notation $AP_n^{q}$ to 
represent the $n$-th $q$-pillow for odd $q$.  Recall from the introduction that this is a centrally symmetric 
region with steps of length $q$ and central belt of size $2\times 2n$.  Calculating the number of tilings of 
$AP_n^5$, $AP_n^7$, and $AP_n^9$ for $n$ at least up to $40$ gives strong evidence that $\#AP_n^{q}$ for $q$ odd 
is always a larger number squared times a smaller number.  As with 3-pillows, we can write 
$\#AP_n^{q}=\ell_{n,q}^2s_{n,q}$, where $s_{n,q}$ is relatively small but is not necessarily the squarefree part 
of $\#AP_n^{q}$.

\begin{figure}
\begin{center}
\begin{tabular}{ccc}
  \epsfig{figure=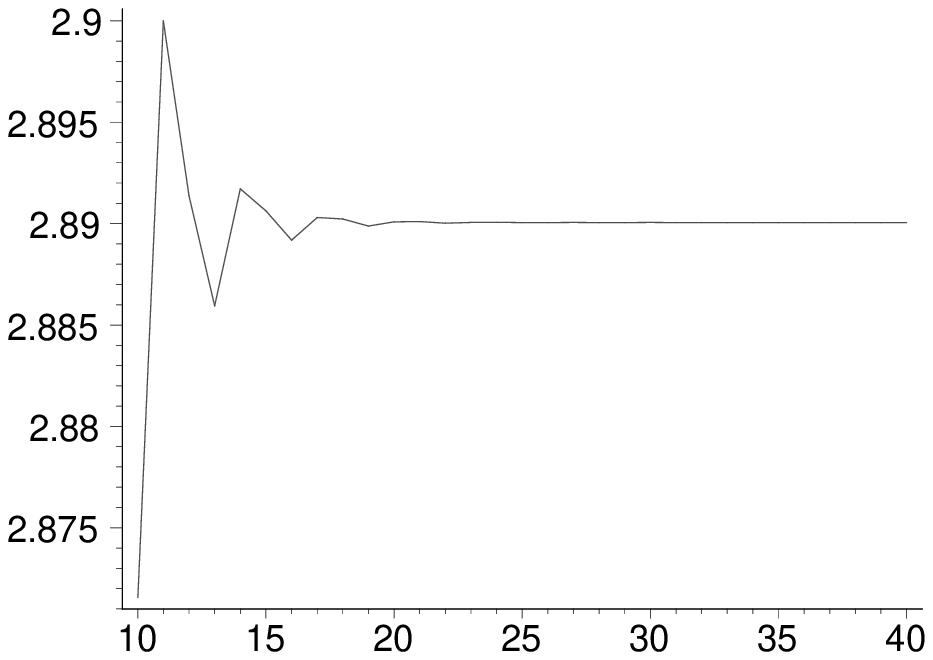, width=2in} &&
  \epsfig{figure=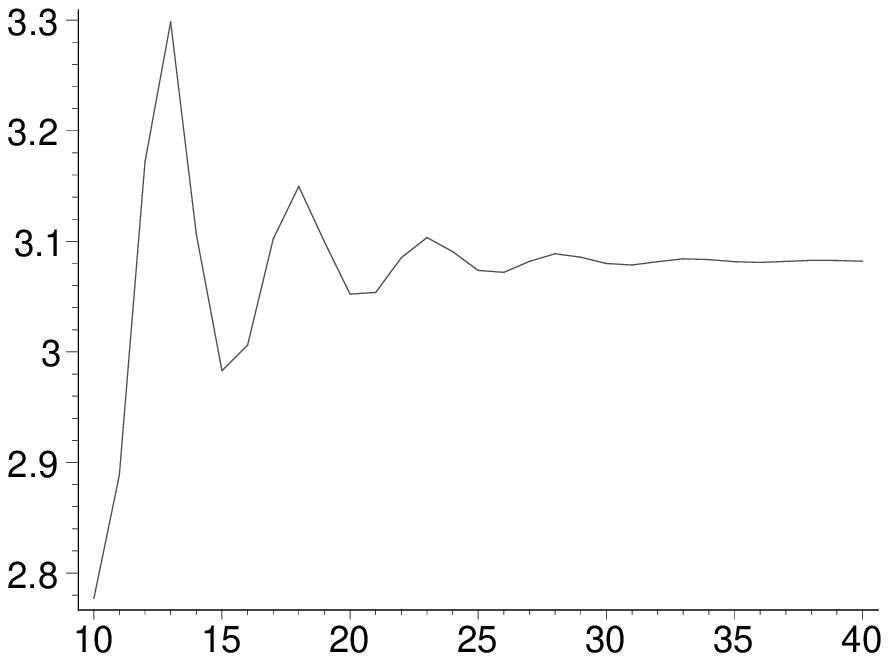, width=2in} \medskip \\
$q=3$ && $q=5$ \\
  \epsfig{figure=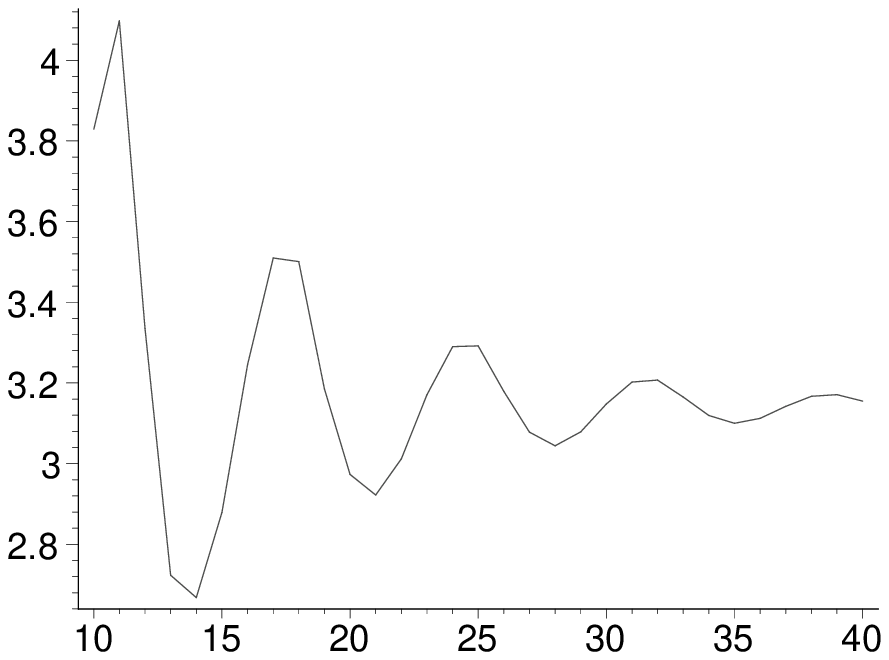, width=2in} &&
  \epsfig{figure=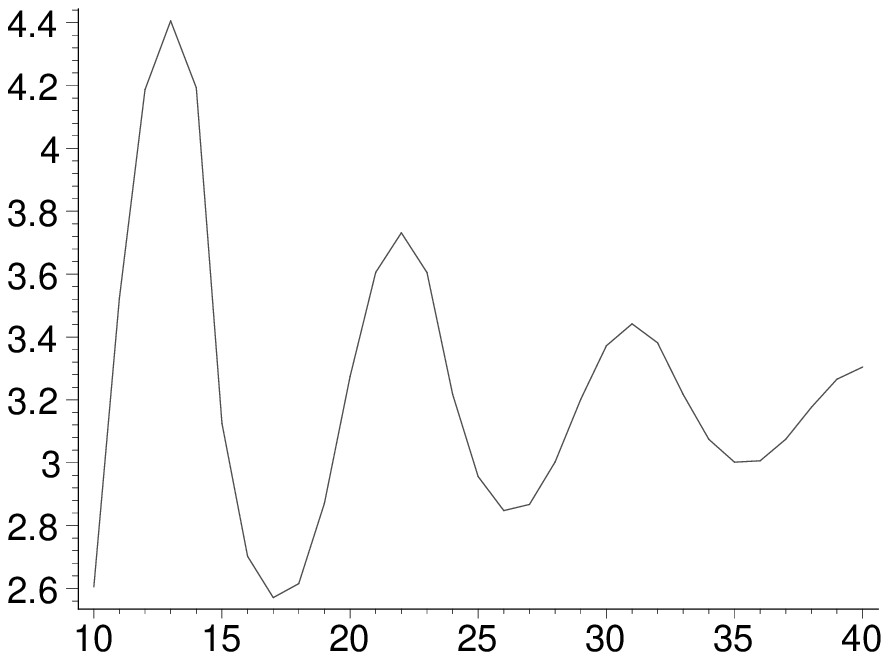, width=2in} \\
$q=7$ && $q=9$ \\
\end{tabular}
\end{center}
\caption{Plots of values of $s_{n,q}/s_{n-2,q}$ for $q=3,5,7,9$}
\label{Plots}
\end{figure}

Additional structure in these smaller numbers also appears.  Propp's 
Conjecture gives an explicit recurrence for the smaller-number values when 
$q=3$.  In the case $q=5$, there is no linear, constant coefficient 
$k$-th degree recurrence for any $k$ up to 20.  However in each case 
($q=5$, $q=7$, and $q=9$) the structure is undeniable.  Plotting the 
values of $s_{n,q}/s_{n-2,q}$ gives a remarkable damped sinusoidal graph 
converging to some number $\rho$.  (See Figure \ref{Plots}.)  

In addition, plotting the values of 
$s_{n,q}/s_{n-1,q}$ depending on whether $n$ is even or odd gives two 
damped sinusoidal graphs converging to some values $\tau$ and $2\tau$ 
respectively.  This implies $\rho=2\tau^2$.  Approximate values for $\rho$ 
are given in Table \ref{RhoValues}.  We have included the values for Aztec 
diamonds ($q=1$) and the original Aztec pillows ($q=3$).

\begin{table}[ht]
\begin{center}
\begin{tabular}{rl}
$q$ & Approximate limit of $s_{n,q}/s_{n-2,q}$ \\ \hline
1 & 2 \\
3 & 2.890053636 \\
5 & 3.0821372 \\
7 & 3.145 \\
9 & 3.18 \\
\end{tabular}
\end{center}
\caption{$\rho$ values for various values of $q$}
\label{RhoValues}
\end{table}

%Appendices \ref{5pillows} and \ref{7pillows} contain the raw data.

\medskip
Even more appears to be true.  By construction, the dual graphs to all odd pillows are 2-even-symmetric---that is, 
centrally symmetric so that antipodes are an even path length apart.  Jockusch \cite{Jockusch} proved that the 
number of matchings of such graphs is the sum of two squares.

Calculating these squares gives us a better understanding of the larger numbers $\ell_n$ appearing in Propp's 
Conjecture, and the larger numbers $\ell_{n,q}$ in general.  In the case of odd pillows, we can write 
$\#AP_n^q=a_{n,q}^2+b_{n,q}^2$ for explicitly known positive integers $a_{n,q}$ and $b_{n,q}$.  Surprisingly, in 
each case I considered, $\ell_{n,q}$ divides both $a_{n,q}$ and $b_{n,q}$.  This implies that the smaller number 
$s_{n,q}$ is also a sum of squares, which may lend some insight into its recurrence formula.

This allows us to formulate a new and improved conjecture about odd pillows.

{\conjecture The number of tilings of $AP_n^{q}$ is a larger number squared times a smaller number.  There exist 
values $\ell_{n,q}$ and $s_{n,q}$ such that we can write $\#AP_n^{q}=\ell_{n,q}^2s_{n,q}$.  The value $s_n$ 
satisfies the following structure.  The ratio of $s_{2n+1,q}/s_{2n,q}$ to $s_{2n+2,q}/s_{2n+1,q}$ is exactly 2 in 
the limit.  In addition, for the explicit values $a_{n,q}$ and $b_{n,q}$ from $\#AP_n^q$'s sum of squares 
representation, the value $\ell_{n,q}$ from above divides both $a_{n,q}$ and $b_{n,q}$. \label{GenConj}}

\subsection{Conclusion}

The Hamburger Theorem allows us to calculate the number of cycle systems in generalized Aztec Pillow graphs more 
quickly. It also allows us to count the number of cycle systems in many previously inaccessible graphs, such as 
non-planar hamburger graphs.  Just as Gessel and Viennot's result has found many applications, I hope the 
Hamburger Theorem will be useful as well.

\subsection{Remarks}

This article is based on work in the author's doctoral dissertation.  A preliminary version of this article 
appeared as a poster at the 2005 International Conference on Formal Power Series and Algebraic Combinatorics in 
Taormina, Italy. The author would like to thank Henry Cohn and Tom Zaslavsky for numerous intriguing discussions 
and useful corrections.

%%%%%%%%%%%%%%%%%%%%%%%%%%%%%%%%%%%%%%%%%%%%%%%%%%%%%%%%%%%%%%%%%%%%%%%%
\nocite{*}   % include everything in the uwthesis.bib file
\bibliographystyle{plain}
\bibliography{2005_11_HamburgerThm}

\providecommand{\noopsort}[1]{}
\begin{thebibliography}{10}

\bibitem{Aigner}
Martin Aigner.
\newblock Lattice paths and determinants.
\newblock In {\em {Computational Discrete Mathematics}}, volume 2122 of {\em
  {Lecture Notes in Comput. Sci.}}, pages 1--12. Springer, Berlin, 2001.

\bibitem{BruKirk}
Richard Brualdi and Stephen Kirkland.
\newblock {Aztec diamonds and digraphs, and Hankel determinants of Schr\"{o}der
  numbers. Submitted, 2003. Available at
  http://www.math.wisc.edu/$\sim$brualdi/aztec2.pdf.}

\bibitem{GV1}
Ira Gessel and Xavier~G. Viennot.
\newblock Binomial determinants, paths, and hook length formulae.
\newblock {\em Adv. in Math.}, 58:300--321, 1985.

\bibitem{GV2}
Ira Gessel and Xavier~G. Viennot.
\newblock Determinants, paths, and plane partitions.
\newblock Manuscript, 1989.
\newblock Available at http://www.cs.brandeis.edu/$\sim$ira/papers/pp.pdf.

\bibitem{Jockusch}
W.~Jockusch.
\newblock Perfect matchings and perfect squares.
\newblock {\em J. Combin. Theory Ser. A}, 67(1):100--115, 1994.

\bibitem{KarMcG}
Samuel Karlin and James McGregor.
\newblock Coincidence probabilities.
\newblock {\em Pacific J. Math.}, 9:1141--1164, 1959.

\bibitem{Kasteleyn}
P.~W. Kasteleyn.
\newblock The statistics of dimers on a lattice~~{I. The} number of dimer
  arrangements on a quadratic lattice.
\newblock {\em Physica}, 27:1209--1225, 1961.

\bibitem{Lindstrom}
B.~Lindstr{\"o}m.
\newblock On the vector representations of induced matroids.
\newblock {\em Bull. London Math. Soc.}, 5:85--90, 1973.

\bibitem{Percus}
J.~Percus.
\newblock One more technique for the dimer problem.
\newblock {\em J. Math. Phys.}, 10:1881--1884, 1969.

\bibitem{Propp}
James Propp.
\newblock Enumeration of matchings: problems and progress.
\newblock In {\em New perspectives in algebraic combinatorics (Berkeley, CA,
  1996--97)}, volume~38 of {\em Math. Sci. Res. Inst. Publ.}, pages 255--291.
  Cambridge Univ. Press, Cambridge, 1999.
\newblock {arXiv:math.CO/9904150}.

\bibitem{EIS}
N.~J.~A. Sloane.
\newblock {The On-Line Encyclopedia of Integer Sequences. \\
  http://www.research.att.com/$\sim$njas/sequences/}.

\bibitem{Zhang}
Fuzhen Zhang.
\newblock {\em {The Schur Complement and its Applications, {\em volume 4 of}
  Numerical Methods and Algorithms}}.
\newblock Springer, New York, 2005.

\end{thebibliography}

\end{document}